\documentclass[11pt,oneside,english]{elsart}
\usepackage[T1]{fontenc}
\usepackage[latin1]{inputenc}
\usepackage{subfigure}
\usepackage{float}
\usepackage{amsmath}
\usepackage{graphicx}
\usepackage{amssymb}

\makeatletter

\providecommand{\tabularnewline}{\\}
\floatstyle{ruled}
\newfloat{algorithm}{tbp}{loa}
\floatname{algorithm}{Algorithm}

\usepackage{color}

\usepackage{algorithm}
\usepackage{algorithmic}

\newcommand{\draftwatermark}{
}

\renewcommand{\draftwatermark}{}

\draftwatermark

\usepackage{babel}
\makeatother
\begin{document}
\begin{frontmatter}

\title{Fast Adaptive Algorithms in the Non-Standard Form for Multidimensional
Problems}

\author{Gregory Beylkin, Vani Cheruvu and Fernando P\'{e}rez }

\address{Department of Applied Mathematics, University of Colorado, Boulder,
CO 80309-0526, United States}

\thanks{This research was partially supported by DOE grant DE-FG02-03ER25583,
DOE/ORNL grant 4000038129 and DARPA/ARO grant W911NF-04-1-0281.}

\begin{abstract}
We present a fast, adaptive multiresolution algorithm for applying
integral operators with a wide class of radially symmetric kernels
in dimensions one, two and three. This algorithm is made efficient
by the use of separated representations of the kernel. We discuss
operators of the class $(-\Delta+\mu^{2}I)^{-\alpha}$, where $\mu\geq0$
and $0<\alpha<3/2$, and illustrate the algorithm for the Poisson
and Schr\"{o}dinger equations in dimension three. The same algorithm
may be used for all operators with radially symmetric kernels approximated
as a weighted sum of Gaussians, making it applicable across multiple
fields by reusing a single implementation.

This fast algorithm provides controllable accuracy at a reasonable
cost, comparable to that of the Fast Multipole Method (FMM). It differs
from the FMM by the type of approximation used to represent kernels
and has an advantage of being easily extendable to higher dimensions. 
\end{abstract}
\end{frontmatter}

\section{Introduction}

For a number of years, the Fast Multipole Method (FMM) \cite{GRE-ROK:1987,CA-GR-RO:1988,CH-GR-RO:1999}
has been the method of choice for applying integral operators to functions
in dimensions $d\le3$. On the other hand, multiresolution algorithms
in wavelet and multiwavelet bases introduced in \cite{BE-CO-RO:1991}
for the same purpose were not efficient enough to be practical in
more than one dimension. Recently, with the introduction of separated
representations \cite{BEY-CRA:2002,BEY-MOH:2002,BEY-MOH:2005}, practical
multiresolution algorithms in higher dimensions \cite{H-F-Y-B:2003,H-F-Y-G-B:2004,Y-F-G-H-B:2004a,Y-F-G-H-B:2004}
became available as well. In this paper we present a new fast, adaptive
algorithm for applying a class of integral operators with radial kernels
in dimensions $d=1,2,3$, and we briefly discuss its extension to
higher dimensions.

In physics, chemistry and other applied fields, many important problems
may be formulated using integral equations, typically involving Green's
functions as their kernels. Often such formulations are preferable
to those via partial differential equations (PDEs). For example, evaluating
the integral expressing the solution of the Poisson equation in free
space (the convolution of the Green's function with the mass or charge
density) avoids issues associated with the high condition number of
a PDE formulation. Integral operators appear in fields as diverse
as electrostatics, quantum chemistry, fluid dynamics and geodesy;
in all such applications fast and accurate methods for evaluating
operators on functions are needed. 

The FMM and our approach both employ approximate representations of
operators to yield fast algorithms. The main difference lies in the
type of approximations that are used. For example, for the Poisson
kernel $1/r$ in dimension $d=3$, the FMM \cite{CH-GR-RO:1999} uses
a plane wave approximation starting from the integral \begin{equation}
\frac{1}{r}=\frac{1}{2\pi}\int_{0}^{\infty}e^{-\lambda(z-z_{0})}\int_{0}^{2\pi}e^{i\lambda((x-x_{0})\cos\alpha+(y-y_{0})\sin\alpha)}d\alpha d\lambda,\end{equation}
where $r=\sqrt{(x-x_{0})^{2}+(y-y_{0})^{2}+(z-z_{0})^{2}}$. The elegant
approximation derived from this integral in \cite{CH-GR-RO:1999}
is valid within a solid angle, and thus requires splitting the application
of an operator into directional regions; the number of such regions
grows exponentially with dimension. For the same kernel, our approach
starts with the integral \begin{equation}
\frac{1}{r}=\frac{2}{\sqrt{\pi}}\int_{-\infty}^{\infty}e^{-r^{2}e^{2s}+s}ds,\label{eq:kernel_integral}\end{equation}
and its discretization with finite accuracy $\epsilon$ yields a spherically
symmetric approximation as a weighted sum of gaussians. Other radial
kernels can be similarly treated by a suitable choice of integrals.
The result is a separated representation of kernels and, therefore,
an immediate reduction in the cost of their application. This difference
in the choice of approximation dictates the differences in the corresponding
algorithms. In dimension $d\le3$ both approaches are practical and
yield comparable performance. The key advantage of our approach is
its straightforward extensibility to higher dimensions \cite{BEY-MOH:2002,BEY-MOH:2005}. 

Given an arbitrary accuracy $\epsilon$, we effectively represent
kernels by a set of exponents and weights describing the terms of
the gaussian approximation of integrals like in (\ref{eq:kernel_integral}).
The number of terms in such sum is roughly proportional to $\log(\epsilon^{-1})$,
or a low power of $\log(\epsilon^{-1})$, depending on the operator.
Since operators are fully described up to an accuracy $\epsilon$
by the exponents and weights of the sum of gaussians, a single algorithm
applies all such operators. These include operators such as $(-\Delta+\mu^{2}I)^{-\alpha}$,
where $\mu\geq0$ and $0<\alpha<3/2$, and certain singular operators
such as the projector on divergence-free functions. Since many physically
significant operators depend only on the distance between interacting
objects, our approach is directly applicable to problems involving
a wide class of operators with radial kernels.

We combine separated and multiresolution representations of kernels
and use multiwavelet bases \cite{ALPERT:1993} that provide \textit{inter
alia} a method for discretizing integral equations, as is the case
in quantum chemistry \cite{H-F-Y-B:2003,H-F-Y-G-B:2004,Y-F-G-H-B:2004,Y-F-G-H-B:2004a}.
This choice of multiresolution bases accommodates integral and differential
operators as well as a wide variety of boundary conditions, without
degrading the order of the method \cite{A-B-G-V:2002,BEYLKI:2001}.
Multiwavelet bases retain the key desired properties of wavelet bases,
such as vanishing moments, orthogonality, and compact support. Due
to the vanishing moments, wide classes of integro-differential operators
have an \textit{effectively sparse} matrix representation, i.e., they
differ from a sparse matrix by an operator with small norm. Some of
the basis functions of multiwavelet bases are discontinuous, similar
to those of the Haar basis and in contrast to wavelets with regularity
(see e.g. \cite{DAUBEC:1992,CHUI:1992}). The usual choices of scaling
functions for multiwavelet bases are either the scaled Legendre or
interpolating polynomials. Since these are also used in the discontinuous
Galerkin and discontinuous spectral elements methods, our approach
may also be seen as an adaptive extension of these methods.

The algorithm for applying an operator to a function starts with computing
its adaptive representation in a multiwavelet basis, resulting in
a $2^{d}$-tree with blocks of coefficients at the leaves. Then the
algorithm adaptively applies the (modified) separated non-standard
form \cite{BE-CO-RO:1991} of the operator to the function by using
only the necessary blocks as dictated by the function's tree representation.
We note that in higher dimensions, $d\gg3$, functions need to be
in a separated representation as well, since the usual constructions
via bases or grids are prohibitive (see \cite{BEY-MOH:2002,BEY-MOH:2005}).

We start in Section~\ref{sec:Preliminary-Considerations} by recalling
the basic notions of multiresolution analysis, non-standard operator
form and adaptive representation of functions underlying our development.
We then consider the separated representation for radially symmetric
kernels in Section~\ref{sec:Separated-Repr-Kernels}, and use it
to efficiently extend the modified ns-form to multiple spatial dimensions
in Section~\ref{sec:Modified-ns-form-Nd}. We pay particular attention
to computing the band structure of the operator based on one dimensional
information. We use this construction in Section~\ref{sec:Multidim-application-nsform}
to introduce the adaptive algorithm for application of multidimensional
operators in the modified ns-form, and illustrate its performance
in Section~\ref{sec:Numerical-examples}. We consider two examples:
the Poisson equation in free space and the ground state of the Hydrogen
atom. We conclude with a brief discussion in Section~\ref{sec:Discussion-and-conclusions}.

\section{\label{sec:Preliminary-Considerations}Preliminary considerations}

This section and Appendix are provided for the convenience of the
reader in order to keep this paper reasonably self-contained. We provide
background material and introduce necessary notation.

The essence of our approach is to decompose the operator using projectors
on a Multiresolution Analysis (MRA), and to efficiently apply its
projections using a separated representation. We use the decomposition
of the operator into the ns-form \cite{BE-CO-RO:1991}, but we organize
it differently (thus, \emph{modified} ns-form) to achieve greater
efficiency. This modification becomes important as we extend this
algorithm to higher dimensions. 

In this section we introduce notation for MRA, describe the adaptive
representation of functions and associated data structures, introduce
the modified ns-form and an algorithm for its adaptive application
in dimension $d=1$ as background material for the multidimensional
case.

\subsection{Multiresolution analysis}

Let us consider the multiresolution analysis as a decomposition of
$L^{2}([0,1]^{d})$ into a chain of subspaces \[
\mathbf{V}_{0}\subset\mathbf{V}_{1}\subset\mathbf{V}_{2}\subset\dots\subset\mathbf{V}_{n}\subset\dots,\]
so that $L^{2}([0,1]^{d})=\overline{\cup_{j=0}^{\infty}\mathbf{V}_{j}}$.
We note that our indexing of subspaces (increasing towards finer scales)
follows that in \cite{A-B-G-V:2002}, and is the reverse of that in
\cite{BE-CO-RO:1991,DAUBEC:1992}. On each subspace $\mathbf{V}_{j}$,
we use the tensor product basis of scaling functions obtained using
the functions $\phi_{kl}^{j}(x)$ $(k=0,\ldots,p-1)$ which we briefly
describe in Appendix.

The wavelet subspaces $\mathbf{W}_{j}$ are defined as the orthogonal
complements of $\mathbf{V}_{j}$ in $\mathbf{V}_{j-1}$, thus \[
\mathbf{V}_{n}=\mathbf{V}_{0}\oplus_{j=0}^{n}\mathbf{W}_{j}\,.\]

Introducing the orthogonal projector on $\mathbf{V}_{j}$, $\mathbf{P}_{j}:L^{2}([0,1]^{d})\to\mathbf{V}_{j}$
and considering an operator $\mathbf{T}:L^{2}([0,1]^{d})\to L^{2}([0,1]^{d})$,
we define its projection $\mathbf{T}_{j}:\mathbf{V}_{j}\to\mathbf{V}_{j}$
as $\mathbf{T}_{j}=\mathbf{P}_{j}\mathbf{T}\mathbf{P}_{j}$. We also
consider the orthogonal projector $\mathbf{Q}_{j}:L^{2}([0,1]^{d})\to\mathbf{W}_{j}$,
defined as $\mathbf{Q}_{j}=\mathbf{P}_{j+1}-\mathbf{P}_{j}$.

\subsection{\label{sec:Adaptive-Repr-Functions}Adaptive representation of functions}

Let us describe an adaptive refinement strategy for construction multiresolution
representations of functions $f:B\rightarrow B$, where $B=[0,1]^{d}$.
We proceed by recursive binary subdivision of the box $B$, so the
basic structure representing our functions is a $2^{d}-$tree with
arrays of coefficients stored at the leaves (terminal nodes) and no
data stored on internal nodes. On each box obtained via this subdivision,
our basis is a tensor product of orthogonal polynomials of degree
$k=0,\ldots,p-1$ in each variable, as described in Appendix~\ref{sec:Appendix}.
Therefore, the leaves carry $d$-dimensional arrays of $p^{d}$ coefficients
which may be used to approximate function values anywhere in the box
corresponding to the spatial region covered by it, via (\ref{eq:func_interp_1d})
or its equivalent for higher values of $d$. For conciseness, we will
often refer to these $d$-dimensional arrays of coefficients stored
at tree nodes as \emph{function blocks}. 

This adaptive function decomposition algorithm is similar to that
used in \cite{ETH-GRE:2001}. Such construction formally works in
any dimension $d$. However, since its complexity scales exponentially
with $d$, its practical use is restricted to fairly low dimension,
e.g. $d\lesssim4$. In higher dimensions, alternate representation
strategies for functions such as \cite{BEY-MOH:2002,BEY-MOH:2005}
should be considered. In high dimensions, these strategies deal with
the exponential growth of complexity by using controlled approximations
that have linear cost in $d$. 

For simplicity, we will describe the procedure for the one-dimensional
case since the extension to dimensions $d=2,\,3$ is straightforward.
Since we can not afford to construct our representation by starting
from a fine scale (especially in $d=3$), we proceed by successive
refinements of an initial coarse sampling. This approach may result
in a situation where the initial sampling is insufficient to resolve
a rapid change in a small volume; however, in practical applications
such situations are rare and may be avoided by an appropriate choice
of the initial sampling scale.

Let $B_{l}^{j}=[2^{-j}l,2^{-j}(l+1)]$, $l=0,\ldots,2^{j}-1,$ represent
a binary subinterval on scale $j$. We denote by $f_{l}^{j}=\left\{ f_{l}^{j}(x_{k})\right\} _{k=0}^{p-1}$
the vector of values of the function $f$ on the Gaussian nodes in
$B_{l}^{j}$. From these values we compute the coefficients $\left\{ s_{kl}^{j}\right\} _{k=0}^{p-1}$
(see (\ref{eq:fvals2coefs}) in Appendix) and interpolate $f(x)$
for any $x\in B_{l}^{j}$ by using (\ref{eq:func_interp_1d}). We
then subdivide $B_{l}^{j}$ into two child intervals, $B_{2l}^{j+1}$
and $B_{2l+1}^{j+1}$, and evaluate the function $f$ on the Gaussian
nodes in $B_{2l}^{j+1}$ and $B_{2l+1}^{j+1}$. We then interpolate
$f$ by using the coefficients $\left\{ s_{kl}^{j}\right\} _{k=0}^{p-1}$
from their parent interval and denote by $\widetilde{f}_{2l}^{j+1}$
and $\widetilde{f}_{2l+1}^{j+1}$ the vectors of interpolated values
on the two subintervals. Now, if for a given tolerance $\epsilon$
either $\left\Vert f_{2l}^{j+1}-\widetilde{f}_{2l}^{j+1}\right\Vert >\epsilon$
or $\left\Vert f_{2l+1}^{j+1}-\widetilde{f}_{2l+1}^{j+1}\right\Vert >\epsilon$,
we repeat the process recursively for both subintervals, $B_{2l}^{j+1}$
and $B_{2l+1}^{j+1}$; otherwise, we keep the coefficients $\left\{ s_{kl}^{j}\right\} _{k=0}^{p-1}$
to represent the function on the entire interval $B_{l}^{j}$. This
interval then becomes a leaf in our tree.

At this stage we use the $\ell^{\infty}$ norm, thus constructing
an approximation $\tilde{f}$ to the original function $f$ such that
$\left\Vert f-\widetilde{f}\right\Vert _{\infty}<\epsilon$, which
immediately implies that $\left\Vert f-\widetilde{f}\right\Vert _{2}<\epsilon$.
This estimate clearly extends to any dimension. Once the approximation
with $\ell^{\infty}$ norm is constructed, the corresponding tree
may be pruned if an application only requires the approximation to
be valid in the $\ell^{2}$ norm. We start this process on the finest
scale and simply remove all blocks whose cumulative contribution is
below $\epsilon$. Other norms, such as $H_{1},$ can be accommodated
by appropriately weighing the error tolerance with a scale-dependent
factor in the initial (coarse to fine) decomposition process.

The complete decomposition algorithm proceeds by following the above
recipe, starting with an initial coarse scale (typically $j=0$) and
continuing recursively until the stopping criterion is met for all
subintervals. In practice, we choose a stopping scale $j_{\textrm{max}}$,
beyond which the algorithm will not attempt to subdivide any further.
Reaching $j_{\textrm{max}}$ means that the function has significant
variations which are not accurately resolved over an interval of width
$2^{-j_{\textrm{max}}}$ using a basis of order $p$. A pseudo-code
listing of this process is presented as Algorithm \ref{alg:Function-decomposition}.

\begin{algorithm}

\caption{\label{alg:Function-decomposition}Adaptive Function Decomposition.}

\begin{algorithmic}

\STATE Start at a coarse scale, typically $j=0$.

\STATE Recursively, for all boxes $b^j$ on scale $j$, proceed as follows:

\STATE Construct the list $C$ of $2^d$ child boxes $b^{j+1}$ on scale $j+1$.

\STATE Compute the values of the function $f(b^j)$ at the $p^d$ Gauss-Legendre
quadrature nodes in $b^j$.

\FORALL{boxes $b^{j+1} \in C$}

  \STATE From $f(b^j)$, interpolate to the Gauss-Legendre quadrature nodes in
  $b^{j+1}$, producing values $\tilde{f}(b^{j+1})$.

  \STATE Compute the values of $f(b^{j+1})$ at the Gauss-Legendre nodes of
  $b^{j+1}$, by direct evaluation.

  \IF{$\left\Vert f(b^{j+1})-\tilde{f}(b^{j+1})\right\Vert_\infty > \epsilon $}

        \STATE Recursively repeat the entire process for all boxes $b^{j+1}
        \in C$.

  \ENDIF

\ENDFOR

\COMMENT{Geting here means that the interpolation from the parent was
successful for all child boxes.  We store the parent's coefficients from $b^j$
in the function tree.}

\end{algorithmic}

\end{algorithm}

\subsubsection{\label{sub:Redundant-tree-structure}Tree structures for representing
functions}

The decomposition Algorithm~\ref{alg:Function-decomposition} naturally
leads to a tree data structure to represent functions, with the leaves
of the tree corresponding to the spatial intervals over which the
multiwavelet basis provides a sufficiently accurate local approximation.
By using (\ref{eq:func_interp_1d}) or its higher-dimensional  extensions,
the only data needed to approximate $f(x)$ anywhere in $B$ is the
array of basis coefficients on these leaves. Thus we use a tree structure
where the leaves store these coefficients and the internal nodes do
not contain any data (and are effectively removed since we use hash
tables for storage). We will refer to this structure as an \emph{adaptive
tree.} Each level in the tree corresponds to a scale in the MRA, with
the root node corresponding to the coarsest projection $f_{0}\in\mathbf{V}_{0}$.

\begin{figure}

\begin{centering}

  \subfigure[Function.]{\includegraphics[width=0.5\textwidth]{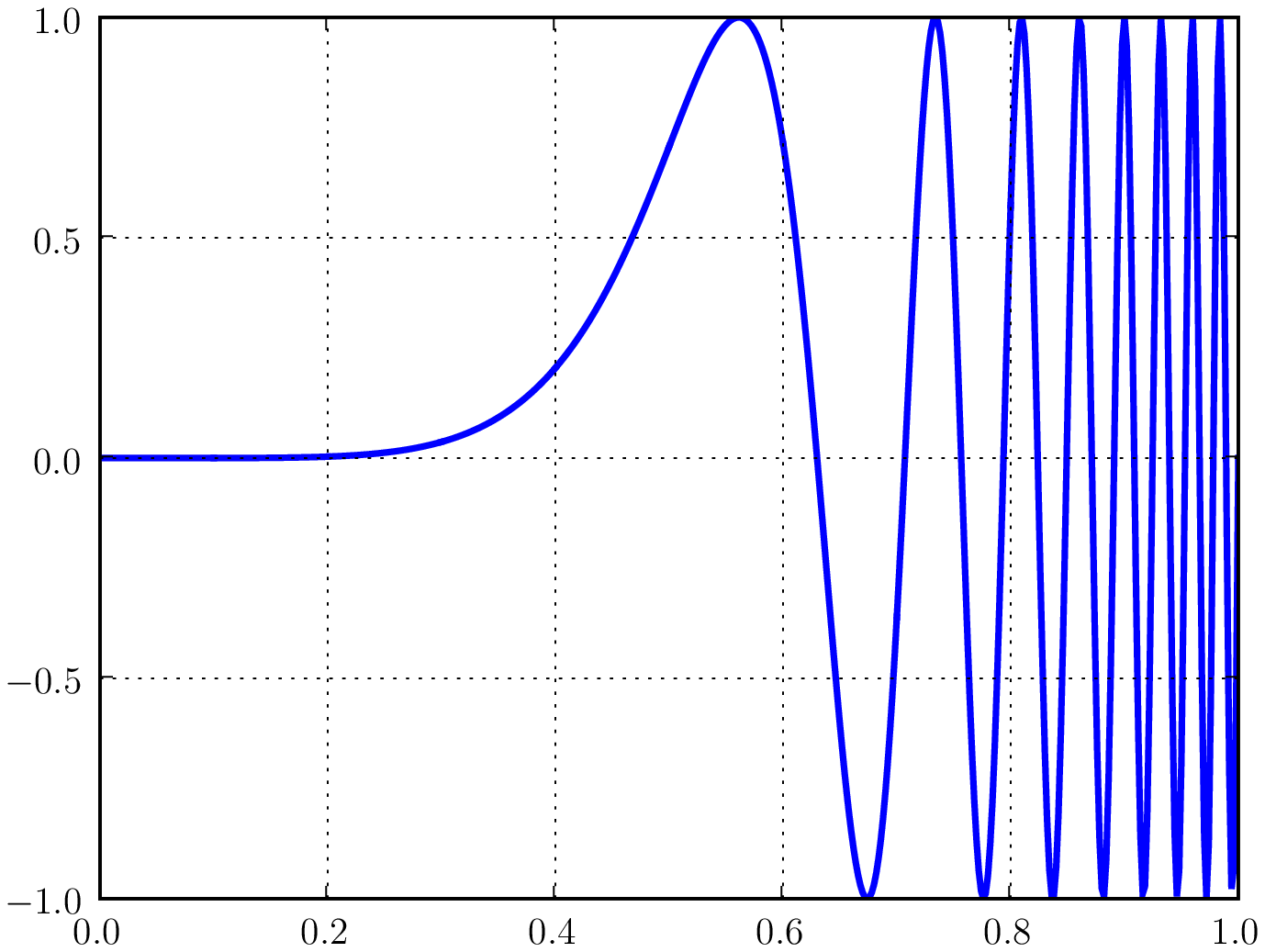}}\subfigure[Pointwise error.]{\includegraphics[width=0.5\textwidth]{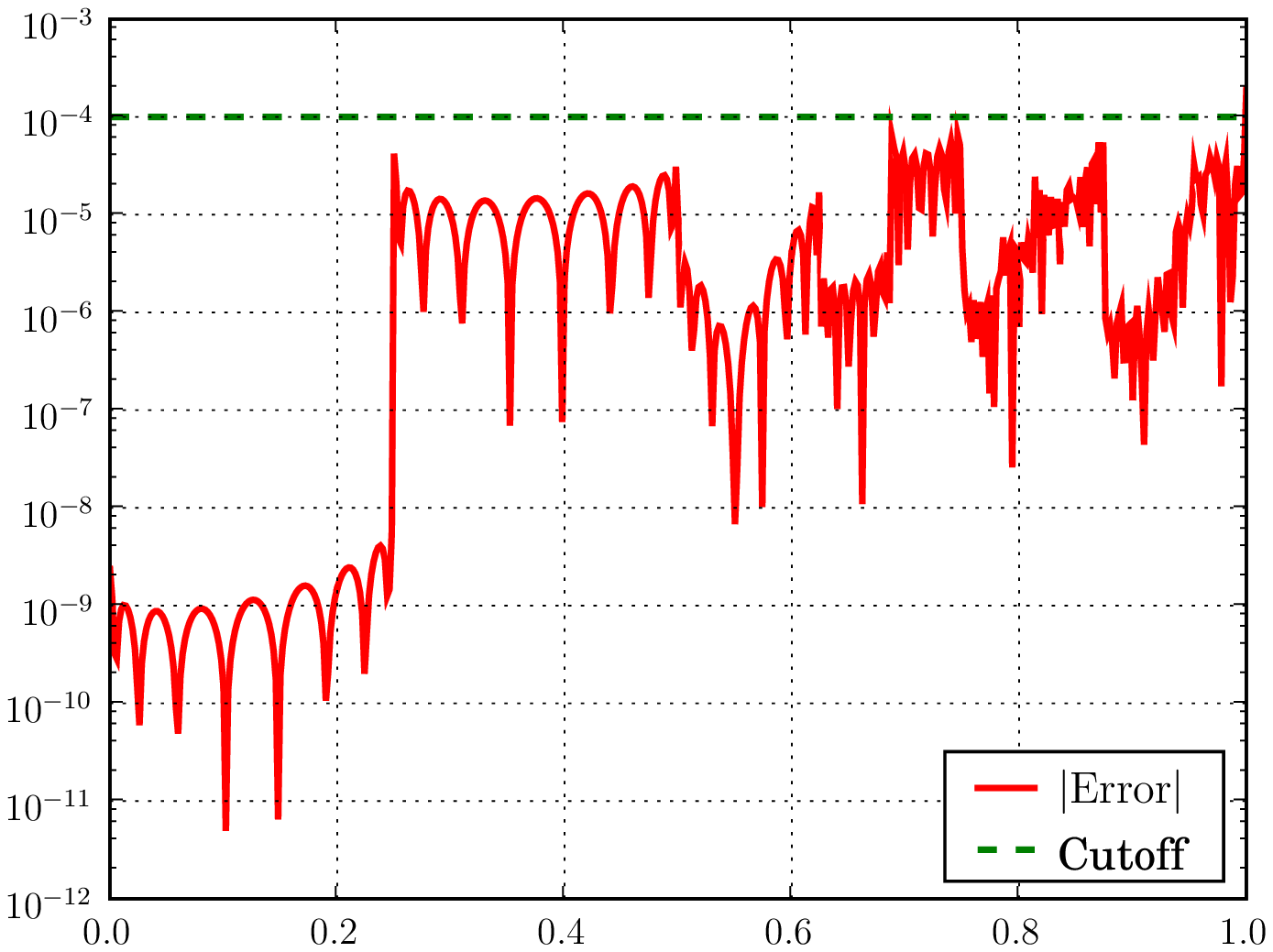}}
\par
\end{centering}

\begin{centering}
\hspace{5mm}
\subfigure[Adaptive tree]{\includegraphics[width=0.4\textwidth]{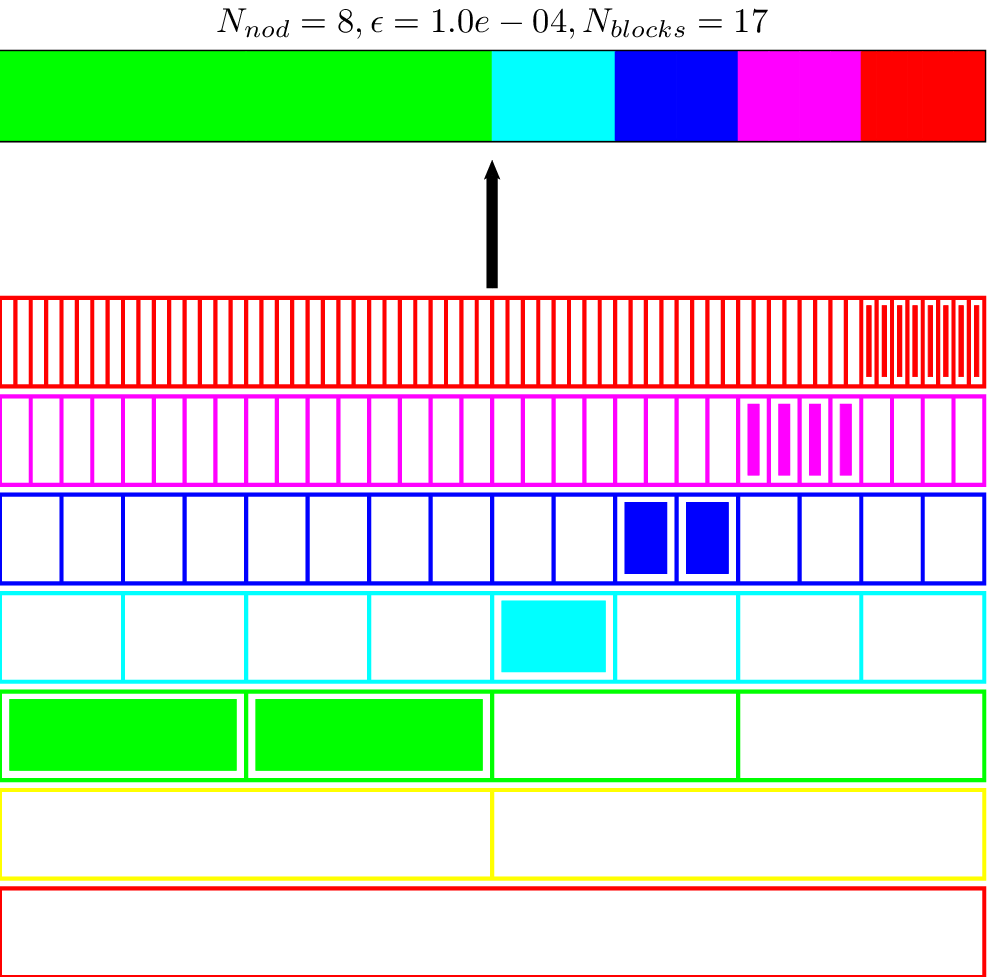}}\hspace{15mm}\subfigure[Redundant tree]{\includegraphics[width=0.4\textwidth]{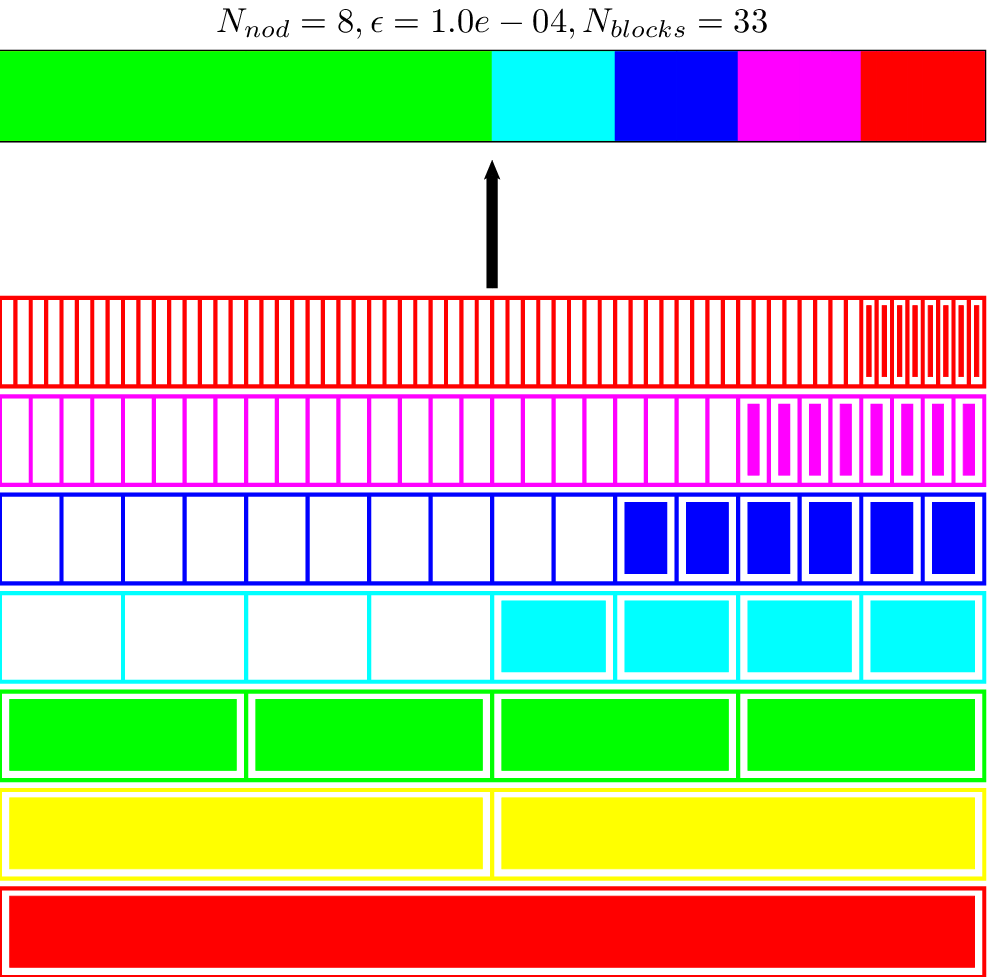}}

\par
\end{centering}

\caption{\label{fig:fun_trees}The function $f(x)=\sin(16\pi x^{6})$ shown
in (a), is decomposed with $p=8$ and $\epsilon=10^{-4}$; (b) shows
the pointwise approximation error. (c) is the resulting adaptive tree,
where smaller subdivisions are required in regions with higher frequency
content. (d) is the redundant tree associated with this adaptive decomposition,
where all internal nodes have been filled with data.}
\end{figure}

Now, in order to apply the modified non-standard form of an operator
to a function, we will show in the next section that we also need
the basis coefficients corresponding to internal nodes of the tree.
Hence, from the adaptive tree data structure we will compute a similar
tree but where we do keep the coefficients of scaling functions on
\emph{all} nodes (leaves and internal).

The coefficients on the internal nodes are redundant since they are
computed from the function blocks stored in the leaves. We will thus
refer to the tree containing coefficients on all nodes as a \emph{redundant
tree}. It is constructed recursively starting from the leaves, by
projecting the scaling coefficients from all sibling nodes onto their
parent node, using the decomposition (\ref{eq:decomp_scaling}).

Figure~\ref{fig:fun_trees} shows both the adaptive and the redundant
trees for a sample function. This figure displays the coarsest scales
at the bottom and progressively finer ones further up, with filled
boxes representing nodes where scaling coefficients are stored and
empty boxes indicating nodes with no data in them (these do not need
to be actually stored in the implementation).

\subsection{\label{sec:Modified-ns-form}Modified ns-form}

The non-standard form \cite{BE-CO-RO:1991} (see also \cite{A-B-G-V:2002}
for the version specialized for multiwavelets) of the operator $\mathbf{T}$
is the collection of components of the telescopic expansion \begin{equation}
\mathbf{T}_{n}=(\mathbf{T}_{n}-\mathbf{T}_{n-1})+(\mathbf{T}_{n-1}-\mathbf{T}_{n-2})+\dots+\mathbf{T}_{0}=\mathbf{T}_{0}+\sum_{j=0}^{n-1}(\mathbf{A}_{j}+\mathbf{B}_{j}+\mathbf{C}_{j}),\label{eq:telescopic}\end{equation}
where $\mathbf{A}_{j}=\mathbf{Q}_{j}\mathbf{T}\mathbf{Q}_{j}$, $\mathbf{B}_{j}=\mathbf{Q}_{j}\mathbf{T}\mathbf{P}_{j}$,
and $\mathbf{C}_{j}=\mathbf{P}_{j}\mathbf{T}\mathbf{Q}_{j}$. The
main property of this expansion is that the rate of decay of the matrix
elements of the operators $\mathbf{A}_{j}$, $\mathbf{B}_{j}$ and
$\mathbf{C}_{j}$ away from the diagonal is controlled by the number
of vanishing moments of the basis and, for a finite but arbitrary
accuracy $\epsilon$, the matrix elements outside a certain band can
be set to zero resulting in an error of the norm less than $\epsilon$.
Such behavior of the matrix elements becomes clear if we observe that
the derivatives of kernels of Calder\'{o}n-Zygmund and pseudo-differential
operators decay faster than the kernel itself. If we use the Taylor
expansion of the kernel to estimate the matrix elements away from
the diagonal, then the size of these elements is controlled by a high
derivative of the kernel since the vanishing moments eliminate the
lower order terms \cite{BE-CO-RO:1991}. We note that for periodic
kernels the band is measured as a periodic distance from the diagonal,
resulting in filled-in `corners' of a matrix representation.

Let us introduce notation to show how the telescopic expansion (\ref{eq:telescopic})
is used when applying an operator to a function. If we apply the projection
of the operator $\mathbf{T}_{j-1}$ not on its ``natural scale'' $j-1,$
but on the finer scale $j$, we denote its upsampled version as $\uparrow(\mathbf{T}_{j-1}).$
In the matrix representation of $\mathbf{T}_{j-1}$, this operation
results in the doubling of the matrix size in each direction. This
upsampling $\uparrow(\cdot)$ and downsampling $\downarrow(\cdot)$
notation will also be used for projections of functions.

With this notation, computing $g=\mathbf{T}f$ via (\ref{eq:telescopic})
splits across scales,\begin{eqnarray}
\hat{g}_{0} & = & \mathbf{T}_{0}f_{0}\nonumber \\
\hat{g}_{1} & = & [\mathbf{T}_{1}-\uparrow(\mathbf{T}_{0})]f_{1}\nonumber \\
\hat{g}_{2} & = & [\mathbf{T}_{2}-\uparrow(\mathbf{T}_{1})]f_{2}\label{eq:g_scales}\\
\ldots & \ldots & \ldots\nonumber \\
\hat{g}_{j} & = & [\mathbf{T}_{j}-\uparrow(\mathbf{T}_{j-1})]f_{j}\nonumber \\
\ldots & \ldots & \ldots\nonumber \end{eqnarray}
where $f_{j}=\mathbf{P}_{j}f$. 

As in the application of the usual ns-form in \cite{BE-CO-RO:1991},
to obtain $g_{n}$ after building the set $\{\hat{g}_{0},\hat{g}_{1},\ldots,\hat{g}_{n}\},$
we have to compute \begin{equation}
g_{n}=\hat{g}_{n}+\uparrow\left(\hat{g}_{n-1}+\uparrow\left(\hat{g}_{n-2}+\uparrow\left(\hat{g}_{n-3}+\ldots+\left(\uparrow\hat{g}_{0}\right)\ldots\right)\right)\right).\label{eq:g_assemble}\end{equation}
The order of the parentheses in this expression is essential, as it
indicates the order of the actual operations which are performed starting
on the coarsest subspace $\mathbf{V}_{0}.$ For example, if the number
of scales $n=4,$ then (\ref{eq:g_assemble}) yields $g_{4}=\hat{g}_{4}+\uparrow\left(\hat{g}_{3}+\uparrow\left(\hat{g}_{2}+\uparrow\left(\hat{g}_{1}+\left(\uparrow\hat{g}_{0}\right)\right)\right)\right)$,
describing the sequence of necessary operations.

Unfortunately, the sparsity of the non-standard form induced by the
vanishing moments of bases is not sufficient for fast practical algorithms
in dimensions other than $d=1$. For algorithms in higher dimensions,
we need an additional structure for the remaining non-zero coefficients
of the representation. We will use separated representations (see
Section~\ref{sec:Separated-Repr-Kernels}) introduced in \cite{BEY-MOH:2002,BEY-MOH:2004P}
and first applied in a multiresolution setting in \cite{H-F-Y-B:2003,H-F-Y-G-B:2004,Y-F-G-H-B:2004a,Y-F-G-H-B:2004}.
Within the retained bands, the components of the non-standard form
are stored and applied in a separated representation and, as a result,
the numerical application of operators becomes efficient in higher
dimensions.

\subsection{Modified ns-form in 1D}

Let us describe a one-dimensional construction for operators on $L^{2}([0,1])$
to introduce all the features necessary for a multidimensional algorithm.
Since we use banded versions of operators, we need to introduce the
necessary bookkeeping. 

The ``template'' for the band structure on scale $j$ comes from the
band on the previous scale $j-1$. For each block on scale $j-1$,
the upsampling operation $\uparrow(\mathbf{T}_{j-1})$ creates four
blocks (all combinations of even/odd row and column indices). We insist
on maintaining the strict correspondence between these four blocks
and those of $\mathbf{T}_{j}$. For this reason the description of
the retained blocks of $\mathbf{T}_{j}$ involves the parity of their
row and column indices. Let us denote the blocks in the matrix representing
$\mathbf{T}_{j}$ by $t^{j;ll'}$, where $l,l'=0\ldots2^{j-1}$. Individual
elements within these blocks are indexed as $t_{ii'}^{j;ll'},$ where
$i,i'=0,\ldots,p-1$, and $p$ is the order of the multiwavelet basis.
For a given width of the band $b_{j},$ we keep the operator blocks
$t^{j;ll'}$ with indices satisfying \begin{equation}
\begin{array}{lrcl}
l-b_{j}+1 & \leq l'\leq & l+b_{j}, & \textrm{for even $l$,}\\
l-b_{j} & \leq l'\leq & l+b_{j}-1, & \textrm{for odd $l$}.\end{array}\label{eq:fine_scale_band}\end{equation}

We denote the banded operators where we keep only blocks satisfying
(\ref{eq:fine_scale_band}) as $\mathbf{T}_{j}^{b_{j}}$ and $\uparrow(\mathbf{T}_{j-1})^{b_{j}}.$
If we downsample the operator $\uparrow(\mathbf{T}_{j-1})^{b_{j}}$
back to its original scale $j-1$, then (\ref{eq:fine_scale_band})
leads to the band described by the condition \begin{equation}
l-\left\lfloor b_{j}/2\right\rfloor \leq l'\leq l+\left\lfloor b_{j}/2\right\rfloor ,\label{eq:down_band}\end{equation}
 where $\left\lfloor b_{j}/2\right\rfloor $ denotes the integer part
of $b_{j}/2.$ We denote the banded operator on scale $j-1$ as $\mathbf{T}_{j-1}^{\left\lfloor b_{j}/2\right\rfloor }$,
where we retain blocks satisfying (\ref{eq:down_band}).

If we now rewrite (\ref{eq:g_scales}) \emph{keeping only blocks within
the bands} on each scale, we obtain\begin{eqnarray}
\hat{g}_{0} & = & \mathbf{T}_{0}f_{0}\nonumber \\
\hat{g}_{1} & = & [\mathbf{T}_{1}^{b_{1}}-\uparrow(\mathbf{T}_{0})^{b_{1}}]f_{1}=\mathbf{T}_{1}^{b_{1}}f_{1}-\uparrow(\mathbf{T}_{0})^{b_{1}}f_{1}\nonumber \\
\hat{g}_{2} & = & [\mathbf{T}_{2}^{b_{2}}-\uparrow(\mathbf{T}_{1})^{b_{2}}]f_{2}=\mathbf{T}_{2}^{b_{2}}f_{2}-\uparrow(\mathbf{T}_{1})^{b_{2}}f_{2}\label{eq:g_scales_band}\\
\ldots & \ldots & \ldots\nonumber \\
\hat{g}_{j} & = & [\mathbf{T}_{j}^{b_{j}}-\uparrow(\mathbf{T}_{j-1})^{b_{j}}]f_{j}=\mathbf{T}_{j}^{b_{j}}f_{j}-\uparrow(\mathbf{T}_{j-1})^{b_{j}}f_{j}\nonumber \\
\ldots & \ldots & \ldots\nonumber \end{eqnarray}

For any arbitrary but finite accuracy, instead of applying the full
$[\mathbf{T}_{j}-\uparrow(\mathbf{T}_{j-1})],$ we will only apply
its banded approximation. 

A simple but important observation is that\begin{equation}
\downarrow([\uparrow(\mathbf{T}_{j-1})]f_{j})=\mathbf{T}_{j-1}f_{j-1}\,,\label{eq:down_up}\end{equation}
which follows from the fact that $\mathbf{Q}_{j}\mathbf{P}_{j}=\mathbf{P}_{j}\mathbf{Q}_{j}=0,$
since these are orthogonal projections. Thus, we observe that $\downarrow(\uparrow(\mathbf{T}_{j-1})^{b_{j}})=\mathbf{T}_{j-1}^{\left\lfloor b_{j}/2\right\rfloor };$
so instead of applying $\uparrow(\mathbf{T}_{j-1})^{b_{j}}f_{j}$
on scale $j$, we can obtain the same result using (\ref{eq:down_up}),
so that $\uparrow\left(\mathbf{T}_{j-1}^{\left\lfloor b_{j}/2\right\rfloor }f_{j-1}\right)=\uparrow(\mathbf{T}_{j-1})^{b_{j}}f_{j}$.
Therefore, we will compute only $\mathbf{T}_{j-1}^{\left\lfloor b_{j}/2\right\rfloor }f_{j-1}$
on scale $j-1$ and combine it with computing $\mathbf{T}_{j-1}^{b_{j-1}}f_{j-1}$.
Incorporating this into (\ref{eq:g_assemble}), we arrive at\begin{eqnarray}
g_{n} & = & \mathbf{T}_{n}^{b_{n}}f_{n}+\uparrow\left[\left(\mathbf{T}_{n-1}^{b_{n-1}}-\mathbf{T}_{n-1}^{\left\lfloor b_{n}/2\right\rfloor }\right)f_{n-1}+\right.\label{eq:apply_algo}\\
 &  & \left.\left.\uparrow\left[\left(\mathbf{T}_{n-2}^{b_{n-2}}-\mathbf{T}_{n-2}^{\left\lfloor b_{n-1}/2\right\rfloor }\right)f_{n-2}+\ldots+\left[\uparrow\left[\left(\mathbf{T}_{0}-\mathbf{T}_{0}^{\left\lfloor b_{1}/2\right\rfloor }\right)f_{0}\right]\right]\ldots\right.\right]\right].\nonumber \end{eqnarray}

Using this expression yields an efficient algorithm for applying an
operator, as on each scale $j$, $\mathbf{T}_{j}^{b_{j}}-\mathbf{T}_{j}^{[b_{j+1}/2]}$
is a sparse object, due to the cancellation which occurs for most
of the blocks. In particular, $\mathbf{T}_{j}^{b_{j}}-\mathbf{T}_{j}^{[b_{j+1}/2]}$
is missing the blocks near the diagonal, and we will refer to it as
an \emph{outer band matrix}. We will call $\mathbf{T}_{j}^{b_{j}}$
a \emph{whole band matrix} as it contains both the inner and outer
bands.

\begin{figure}
\begin{centering}

\includegraphics[width=0.8\textwidth]{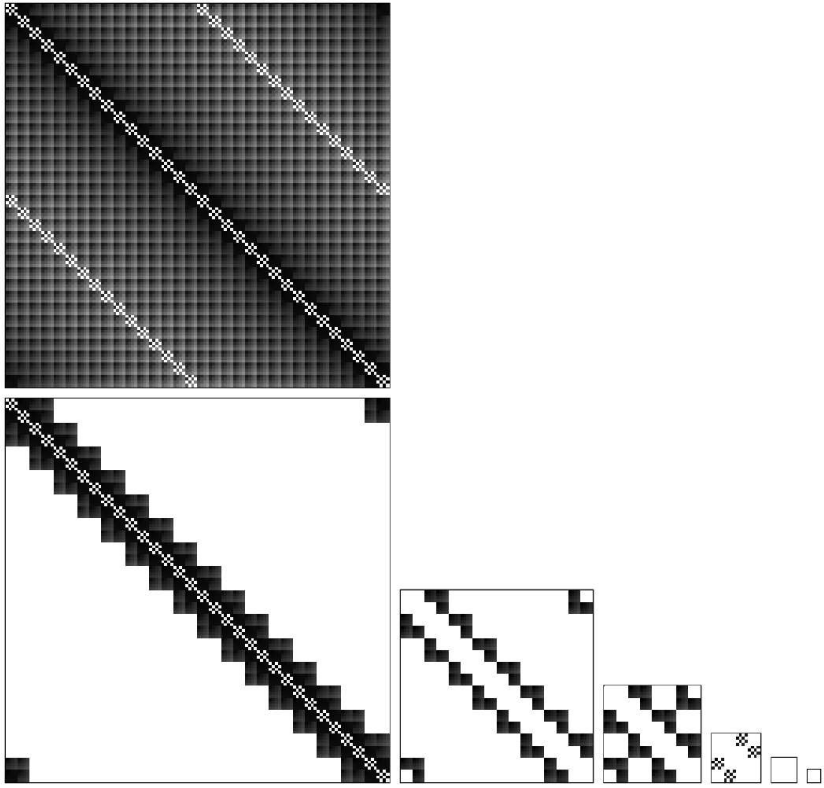}

\par\end{centering}

\caption{\label{fig:oper_bands_1d}Modified non-standard form of the convolution
operator in (\ref{CotangetKernel}) in a multiwavelet basis, with
white representing 0 and black representing large values. The top
matrix is the projection of this operator on $\mathbf{V}_{5}$, resulting
in a dense matrix. The lower half depicts the multiresolution representation
in (\ref{eq:apply_algo}) with only the blocks that are actually retained
for a given accuracy. We will call the leftmost matrix in this series
a \emph{whole band matrix} and all others \emph{outer band matrices.}
The two empty outer band matrices on scales $j=0,1$ are explained
in the main text.}
\end{figure}

The structure of these matrices is illustrated in Figure~\ref{fig:oper_bands_1d}.
The two empty scales $j=0,1$ arise due to the complete cancellation
of blocks on these scales. We note that the modified non-standard
form is constructed adaptively in the number of scales necessary for
a given function. For just two scales, this construction will have
the scale $j=1$ non-empty. Given an adaptive decomposition of a function
on $n$ scales, we precompute the modified non-standard form (depicted
in Figure~\ref{fig:oper_bands_1d} for five scales) on all scales
$n,n-1,\dots,1$. For matrices requiring $2^{n}\cdot2^{n}$ blocks
on the finest scale $n$, we need to keep and apply only $\mathcal{O}(2^{n})$
blocks, as with the original non-standard form in \cite{BE-CO-RO:1991}.

\subsubsection{\label{sec:Adaptive-Application-1d}Adaptive application}

Let us show how to use the multiscale representation in (\ref{eq:apply_algo})
to apply the operator $T$ to a function $f$ with controlled accuracy
$\epsilon$. We describe an adaptive application of the operator to
a function, where we assume (as is often the case) that the tree structure
of the input is sufficient to adequately describe the output with
accuracy $\epsilon$. This assumption will be removed later.

Our algorithm uses the structure shown in Figure~\ref{fig:oper_bands_1d}
in an adaptive fashion. We copy the structure of the redundant tree
for the input function, and use that as a template to be filled for
the output $g$. For each node of the output tree we determine whether
it is a leaf or an internal node: for leaves, we must apply a whole
band \emph{}matrix on the scale of that node, such as the leftmost
matrix for $j=5$ in the example shown in Figure~\ref{fig:oper_bands_1d}.
For internal nodes, we apply an outer band matrix (for that scale).
We note that our construction of the operator produces both whole
and outer band matrices for all scales, and we simply choose the appropriate
kind for each node of output as needed. Upon completion of this process,
we apply the projection (\ref{eq:g_assemble}) to construct the final
adaptive tree representing the output.

\begin{algorithm}

\caption{\label{alg:op-apply-1d}Adaptive non-standard form operator application
in $d=1$, $g=\mathbf{T}f$}

\begin{algorithmic}

\STATE {\bf Initialization:} Construct the redundant tree for $f$ and copy it
as skeleton tree for $g$ (see Section~\ref{sub:Redundant-tree-structure}).

\FORALL{scales $j=0,\ldots,n-1$}

  \FORALL{function blocks $g_l^j$ in the tree for $g$ at scale $j$}

  \STATE \COMMENT{Step 1. Determine the list of all contributing blocks of the
  modified ns-form ${\mathbf{T}_{ll'}^{j}}$ (see
  Section~\ref{sec:Modified-ns-form}):}

  \IF{$g_l^j$ belongs to a \emph{leaf}}

    \STATE Read operator blocks ${\mathbf{T}_{ll'}^{j}}$ from row $l$ of
    \emph{whole band matrix} $\mathbf{T}_{j}^{b_{j}}$.

  \ELSE

    \STATE Read operator blocks ${\mathbf{T}_{ll'}^{j}}$ from row $l$ of
    \emph{outer band matrix} $\mathbf{T}_{j}^{b_{j}}- \mathbf{T}_{j}^{[b_{j+1}/2]}$.

  \ENDIF

  \STATE \COMMENT{Step 2. Find the required blocks ${f_{l'}^j}$ of the input
    function $f$:}

  \IF {function block ${f_{l'}^j}$ exists in the redundant tree for $f$}
    \STATE Retrieve it.
  \ELSE
    \STATE Create it by interpolating from a coarser scale and cache for
reuse.
  \ENDIF

  \STATE \COMMENT{Step 3. Output function block computation:}

  \STATE Compute the resulting output function block according to
  $\hat{g}_{l}^{j}=\sum_{l'}\mathbf{T}_{ll'}^{j}f_{l'}^{j}$, where the
  operation $\mathbf{T}_{ll'}^{j}f_{l'}^{j}$ indicates a regular matrix-vector
  multiplication.

\ENDFOR

\ENDFOR

\STATE \COMMENT{Step 4. Adaptive projection:}

\STATE Project resulting output function blocks $\hat{g}_{l}^{j}$ on all scales
into a proper adaptive tree by using Eq.~(\ref{eq:g_assemble}).

\STATE Discard from the resulting tree unnecessary function blocks at the
requested accuracy.

\STATE {\bf Return:} the function $g$ represented by its adaptive tree.

\end{algorithmic}

\end{algorithm}

Algorithm~\ref{alg:op-apply-1d} returns an adaptive tree representing
the function $g$. This tree contains sufficient information to evaluate
$g$ at arbitrary points by interpolation and may be used as an input
in further computations.

We note that Step~2 in Algorithm~\ref{alg:op-apply-1d} naturally
resolves the problem that is usually addressed by mortar methods,
see e.g.~\cite{MA-MA-PA:1989,A-M-M-P:1990,BE-MA-PA:1993,BE-MA-PA:1994}.
Since adaptive representations have neighboring blocks of different
sizes, they encounter difficulties when applying non-diagonal operators,
as they require blocks which do not exist on that scale. Our approach
simply constructs these as needed and caches them for reuse, without
requiring any additional consideration on the part of the user.

\subsubsection{Numerical example}

\begin{figure}

\begin{centering}

  \includegraphics[width=0.5\textwidth]{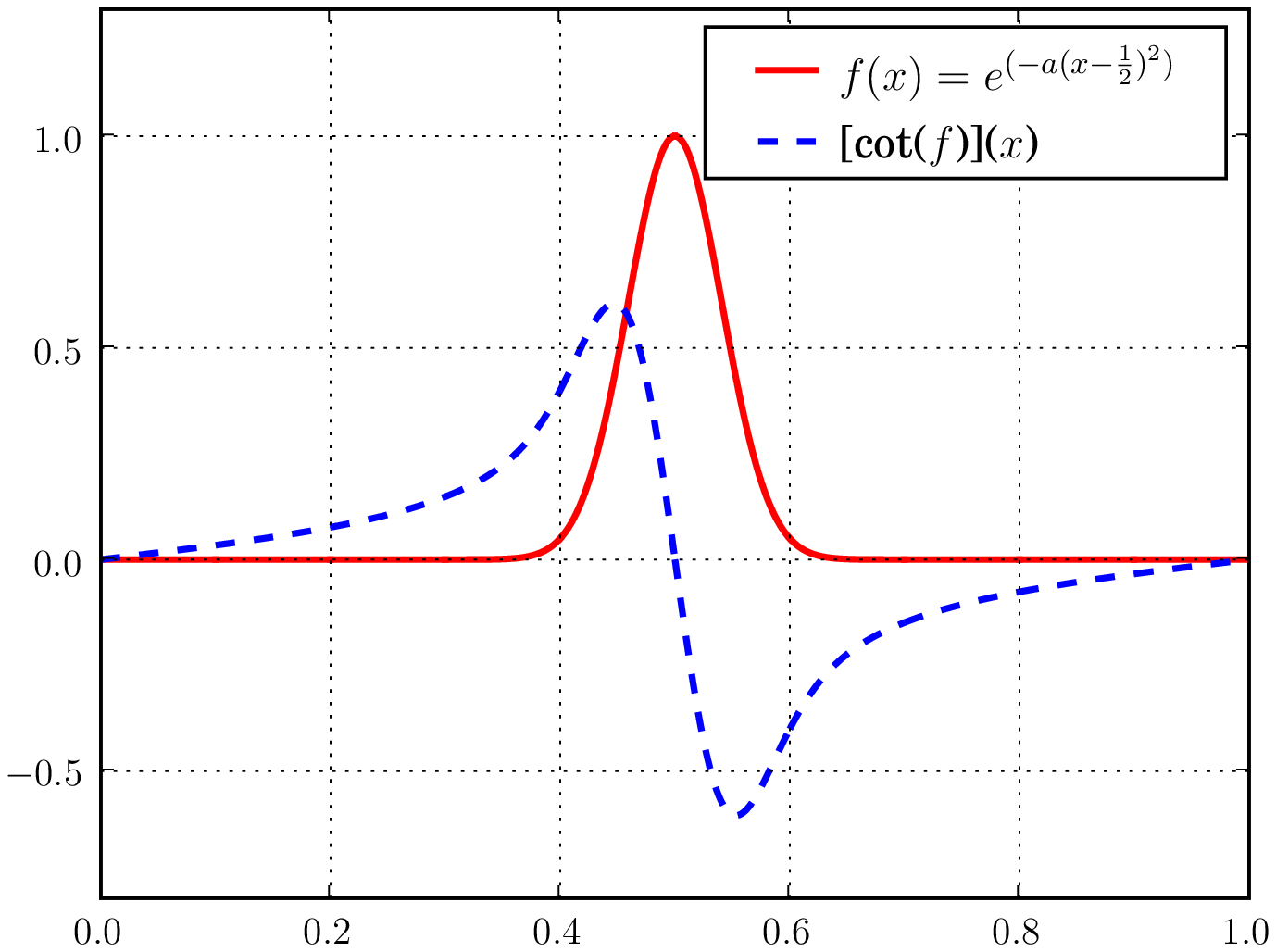}%
  \includegraphics[width=0.5\textwidth]{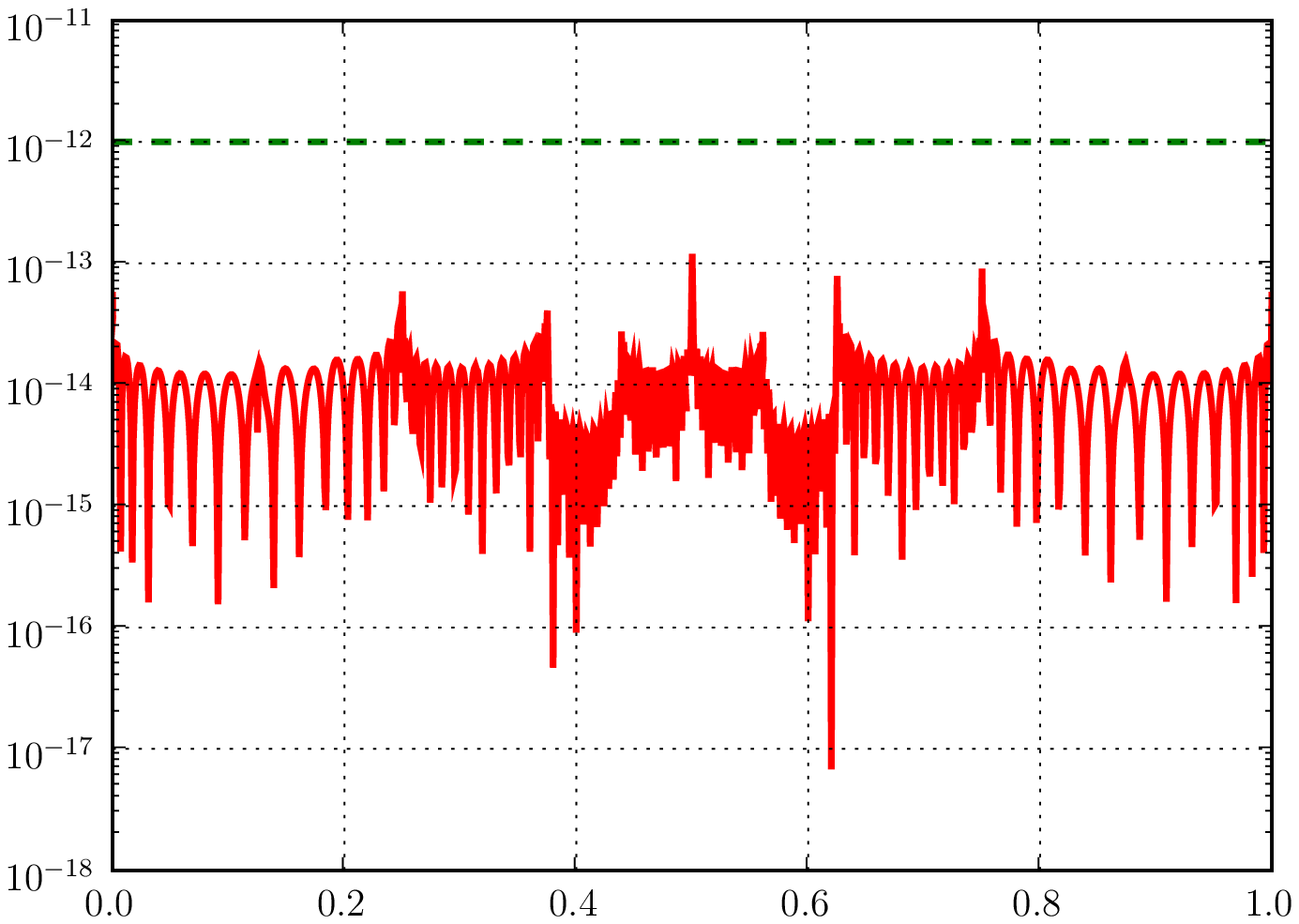}%
  \par
\end{centering}

\caption{\label{fig:cotangent_apply}Results of applying the cotangent kernel
to a periodized Gaussian using basis of order $p=14$ (the last row
in Table~\ref{table:cotangent_apply}). The pointwise error is shown
on the right for a requested accuracy of $\epsilon=10^{-12}$.}
\end{figure}

Let us briefly illustrate the application of the modified ns-form
with an example of a singular convolution on the unit circle, the
operator with the kernel $K(x)=\cot(\pi x),$\begin{equation}
(Cf)(y)\,=\mbox{p.v.}\,\int_{0}^{1}\cot(\pi(y-x))\, f(x)\, dx,\label{CotangetKernel}\end{equation}
a periodic analogue of the Hilbert transform. In order to find its
representation in multiwavelet bases, we compute 

\begin{equation}
r_{{ii'}}^{j;\, l}\,=\,2^{-j}\,\int_{-1}^{1}\, K(2^{-j}(x+l))\,\Phi_{ii'}(x)\, dx\,=2^{-j}\,\int_{-1}^{1}\,\cot(\pi\,2^{-j}(x+l))\,\Phi_{ii'}(x)\, dx\,,\label{cot.01}\end{equation}
where $\Phi_{ii'}(x)$, $i,i'=0,\ldots,k-1$ are cross-correlation
functions described in Appendix~\ref{sub:Cross-correlation-functions-of}
and $l=0,\pm1,\pm2,\ldots2^{j}-1$. We compute $r_{{ii'}}^{j;\, l}$
using the convergent integrals\[
r_{{ii'}}^{j;\, l}\,=\,2^{-j}\,\sum_{k=i'-i}^{i'+i}\, c_{ii'}^{k}\,\int_{0}^{1}\,\Phi_{k,0}^{+}(x)\,\left(\cot(\pi\,2^{-j}(x+l))+(-1)^{i+i'}\cot(\pi\,2^{-j}(-x+l))\right)\, dx,\]
where $\Phi_{k,0}^{+}$ is a polynomial described in Appendix~\ref{sub:Cross-correlation-functions-of}.
In our numerical experiment, we apply (\ref{CotangetKernel}) to the
periodic function on $[0,1]$, \[
f(x)=\sum_{k\in\mathbb{Z}}e^{-a(x+k-1/2)^{2}},\]
 which yields\begin{equation}
(Cf)(y)=-\sum_{k\in\mathbb{Z}}e^{-a(y+k-1/2)^{2}}\mbox{Erfi}[\sqrt{a}(y+k-1/2)]=i\sqrt{\frac{\pi}{a}}\sum_{n\in\mathbb{Z}}\mbox{sign}(n)e^{-n^{2}\pi^{2}/a}e^{2\pi iny},\label{eq:cot_gauss}\end{equation}
where $e^{-y^{2}}\mbox{Erfi}(y)=\frac{2}{\sqrt{\pi}}\int_{0}^{y}e^{s^{2}-y^{2}}ds.$
Expression (\ref{eq:cot_gauss}) is obtained by first observing that
the Hilbert transform of $e^{-ax^{2}}$ is $-e^{-ay^{2}}\mbox{Erfi}(\sqrt{a}y)$,
and then evaluating the lattice sum, noting that (see \cite[formula 4.3.91]{ABR-STE:1970})\[
\cot(\pi x)=\frac{1}{\pi}(\frac{1}{x}+\sum_{k=1}^{\infty}\frac{2x}{x^{2}-k^{2}}).\]
Table~\ref{table:cotangent_apply} summarizes the numerical construction
of this solution for $a=300$, at various requested precisions. Optimal
performance is obtained by adjusting the order of the basis $p$ as
a function of the requested precision, to ensure that the operator
remains a banded matrix with small band, and that the adaptive representation
of the input function requires a moderate number of scales. The resulting
numerical error (as compared to the exact analytical solution), measured
in the $\ell^{2}$ norm, is shown in the last column. Figure~\ref{fig:cotangent_apply}
shows the input and results for this example, as well as the point-wise
error for the case where $\epsilon_{\textrm{req }}=10^{-12}$ and
$p=14$ (the last row in the table). 

\begin{table}
\begin{centering}\begin{tabular}{|c|c|c|c|c|}
\hline 
$p$&
Scales&
$N_{\textrm{blocks}}$&
$\epsilon$&
$E_{2}$\tabularnewline
\hline
\hline 
5&
{[}2,3,4]&
8&
$10^{-3}$&
$1.5\cdot10^{-4}$\tabularnewline
\hline 
8&
{[}2,4,5]&
12&
$10^{-6}$&
$1.3\cdot10^{-7}$\tabularnewline
\hline 
11&
{[}2,4,5]&
14&
$10^{-9}$&
$1.1\cdot10^{-10}$\tabularnewline
\hline 
14&
{[}3,4,5]&
16&
$10^{-12}$&
$4.4\cdot10^{-13}$\tabularnewline
\hline
\end{tabular}\par\end{centering}

\bigskip{}

\caption{\label{table:cotangent_apply}Results from evaluating (\ref{eq:cot_gauss})
with our algorithm. The order of the basis $p$ is adjusted as a function
of the requested precision $\epsilon$. The second column indicates
scales present in the adaptive tree for the input. The third column
shows the total number of blocks of coefficients in this tree. The
last column ($E_{2}$) shows the actual error of the computed solution
in the $\ell^{2}$ norm.}
\end{table}

\section{\label{sec:Separated-Repr-Kernels}Separated representations of integral
kernels}

The approach we've discussed so far does not efficiently generalize
to the application of non-separable multidimensional integral kernels.
Since several physically important kernels belong to this category
(e.g. the Poisson kernels in $d=2$ and $d$=3), additional tools
are needed. We now describe the key idea that allows us to perform
this generalization to $d>1$.

We use the separated representation of operators introduced in \cite{BEY-MOH:2002,BEY-MON:2005}
to reduce the computational cost of the straightforward generalization
of the multiresolution approach in \cite{BE-CO-RO:1991}. Such representations
are particularly simple for convolution operators and are based on
approximating kernels by a sum of Gaussians \cite{BEY-MON:2005,H-F-Y-B:2003,H-F-Y-G-B:2004,Y-F-G-H-B:2004,Y-F-G-H-B:2004a}.
This approximation has a multiresolution character by itself and requires
a remarkably small number of terms. In fact, our algorithm uses the
coefficients and the exponents of the Gaussian terms as the only input
from which it selects the necessary terms, scale by scale, according
the desired accuracy threshold $\epsilon$. Therefore, our algorithm
works for all operators with kernels that admit approximation by a
sum of Gaussians. Examples of such operators include the Poisson and
the bound state Helmholtz Green's functions, the projector on divergence
free functions, as well as regular and fractional derivative operators.
Let us consider a particular family of operators $(-\Delta+\mu^{2}I)^{-\alpha}$,
where $\mu\geq0$ and $0<\alpha<3/2$. The kernel of this operator
\[
K_{\mu,\alpha}(r)=2^{-\frac{3}{2}+\alpha}\cdot C_{\alpha}\cdot(\frac{\mu}{r})^{\frac{3}{2}-\alpha}\textrm{K}_{\frac{3}{2}-\alpha}(\mu r),\]
where $\textrm{K}_{\frac{3}{2}-\alpha}$ is the modified Bessel function,
$r=||\mathbf{x}-\mathbf{y}||$ and $C_{\alpha}=2\cdot\,2^{-2\alpha}\pi^{-\frac{3}{2}}/\Gamma(\alpha)$,
has an integral representation\begin{equation}
K_{\mu,\alpha}(r)=C_{\alpha}\,\int_{-\infty}^{\infty}e^{-r^{2}e^{2s}}e^{-\frac{1}{4}\mu^{2}e^{-2s}+(3-2\alpha)s}ds.\label{IntegralRepr}\end{equation}
Using the trapezoidal rule, we construct an approximation valid over
a range of values $\delta\le r\le R$ with accuracy $\epsilon$, of
the form\begin{equation}
\left|K_{\mu,\alpha}(r)-\sum_{m=1}^{M}w_{m}e^{-\tau_{m}r^{2}}\right|\le\epsilon K_{0,\alpha}(r)=\epsilon\frac{\Gamma(3/2-\alpha)\cdot C_{\alpha}}{2r^{3-2\alpha}},\label{KernelApproxByGaussians}\end{equation}
where $\tau_{m}=e^{2s_{m}}$, $w_{m}=h\, C_{\alpha}\, e^{-\mu^{2}e^{-2s_{m}}/4+(3-2\alpha)s_{m}}$,
$h=(B-A)/M$ and $s_{m}=A+mh$. The limits of integration, $A$, $B$
and the step size $h$ are selected as indicated in \cite{BEY-MON:2005},
where it is shown that for a fixed accuracy $\epsilon$ the number
of terms $M$ in (\ref{KernelApproxByGaussians}) is proportional
to $\log(R\delta^{-1})$. Although it is possible to select $\delta$
and $R$ following the estimates in \cite{B-C-F-H:2007} and optimize
the number of terms for a desired accuracy $\epsilon$, in this paper
we start with an approximation that has an obviously excessive range
of validity and thus, an excessive number of terms. 

\begin{figure}

\begin{centering}\vspace{3mm}
  
  \includegraphics[width=0.6\textwidth]{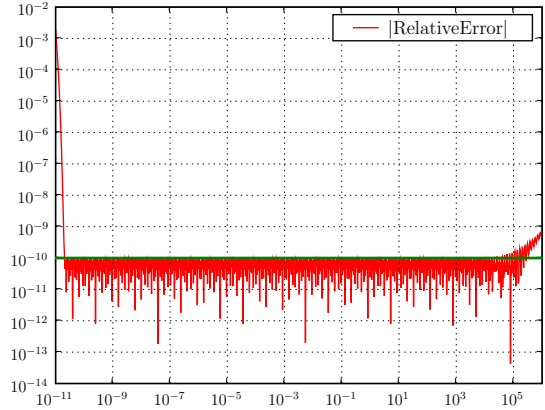}

\par
\end{centering}

\caption{\label{fig:poisson_kernel_error}Relative error of the Gaussian approximation
for the Poisson kernel in 3 dimensions. This unoptimized expansion
uses 299 terms to cover a dynamic range of roughly 15 decades with
$\epsilon=10^{-10}$ relative accuracy.}
\end{figure}

An example of such approximation is shown in Figure~\ref{fig:poisson_kernel_error}.
For a requested tolerance of $\epsilon=10^{-10}$, roughly 300 terms
are enough to provide a valid approximation over a range of $15$
decades. We then let the algorithm choose the necessary terms, scale
by scale, to satisfy the user-supplied accuracy requirement $\epsilon$.
This approach may end up with a few extra terms on some scales in
comparison with that using a nearly optimal number of terms \cite{H-F-Y-B:2003,Y-F-G-H-B:2004a,Y-F-G-H-B:2004,H-F-Y-G-B:2004}.
Whereas the cost of applying a few extra terms is negligible, we gain
significantly in having a much more flexible and general algorithm.

We note that approximation in (\ref{KernelApproxByGaussians}) clearly
reduces the problem of applying the operator to that of applying a
sequence of  Gauss transforms \cite{GRE-STR:1991,GRE-SUN:1998}, one
by one. From this point of view, the algorithm that we present may
be considered as a procedure for applying a linear combination of
Gauss transforms simultaneously.

In order to represent the kernel $K$ of the operator in multiwavelet
bases, we need to compute the integrals, 

\begin{equation}
r_{i_{1}^{\,}i_{1}',i_{2}^{\,}i_{2}',i_{3}^{\,}i_{3}'}^{j;\,{\bf \ell}}\,=\sum_{m=1}^{M}\,2^{-3j}\,\int_{B}\, w_{m}\, e^{-\tau_{m}\Vert2^{-j}({\bf x}+{\bf \ell})\Vert^{2}}\,\Phi_{i_{1}^{\,}i_{1}'}(x_{1})\,\Phi_{i_{2}^{\,}i_{2}'}(x_{2})\,\Phi_{i_{3}^{\,}i_{3}'}(x_{3})\, d{\bf x}\,,\label{KernelApproxByGaussians1}\end{equation}
where $\Phi_{ii'}(x)$ are the cross-correlations of the scaling functions
(see Appendix). We obtain 

\begin{equation}
r_{i_{1}^{\,}i_{1}',i_{2}^{\,}i_{2}',i_{3}^{\,}i_{3}'}^{j;\,{\bf \ell}}\,=\,\sum_{m=1}^{M}\, w_{m}\, F_{i_{1}^{\,}i_{1}'}^{j;\, m,l_{1}}\, F_{i_{2}^{\,}i_{2}'}^{j;\, m,l_{2}}\, F_{i_{3}^{\,}i_{3}'}^{j;\, m,l_{3}}\,,\label{KernelApproxByGaussians2}\end{equation}
where\begin{equation}
F_{ii'}^{j;\, m,l}=\frac{1}{2^{j}}\int_{-1}^{1}e^{-\tau_{m}(x+l)^{2}/4^{j}}\Phi_{ii'}(x)dx.\label{KernelFactor}\end{equation}

Since the Gaussian kernel is not homogeneous, we have to compute integrals
(\ref{KernelFactor}) for each scale. Although in principle $l\in\mathbb{Z}$,
in the next section we explain how to restrict it to a limited range
on each scale, for a given accuracy $\epsilon$.

\section{\label{sec:Modified-ns-form-Nd}Modified ns-form of a multidimensional
operator}

\begin{figure}

\begin{centering}
  \subfigure[Norms of each one-dimensional block computed via
  (\ref{KernelFactor}) for scale $j=1$, as a function of the index $m$ in the 
  separated representation. ]{\includegraphics[width=0.46\textwidth]{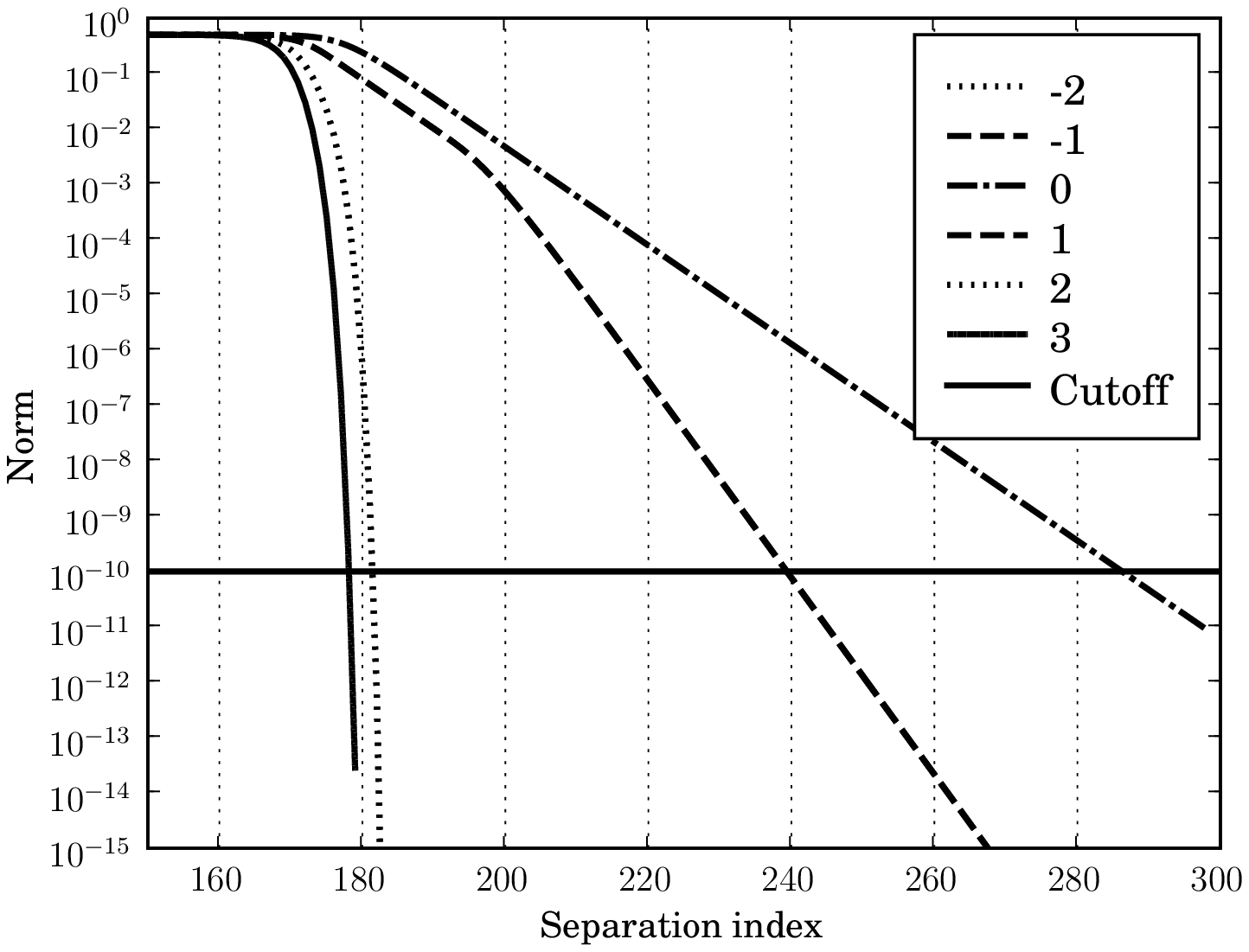}
  }\hspace{8mm}\subfigure[Norm estimates (\ref{NormEstimateUsing1DBlocks}) for
  scale $j=1$ as a function of the index $m$ in the separated
  representation. Based on these estimates, only terms above the cutoff are
  actually applied.] {\includegraphics[width=0.46\textwidth]{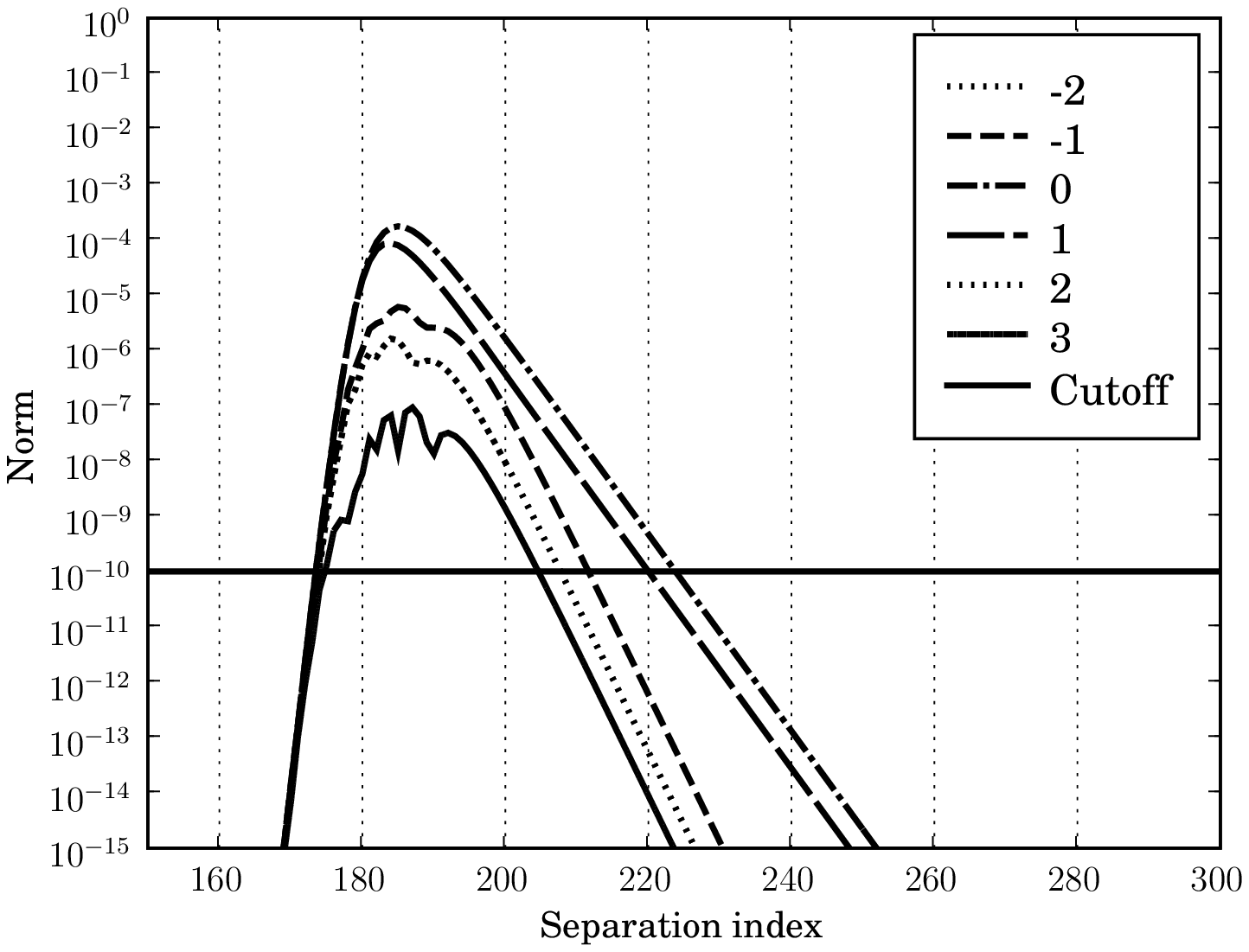} }\par
\end{centering}

\caption{\label{fig:poisson_kernel_norms}Comparison of norms of matrix blocks
generated by the Gaussian approximation for the Poisson kernel in
dimension $d=3$. In each picture, the curves correspond to the different
offsets $l$ for which blocks are generated. Figure (b) illustrates
the estimate in (\ref{NormEstimateUsing1DBlocks}) (see main text
for details).}
\end{figure}

In this section, we describe how the separated representation approximations
of Section~\ref{sec:Separated-Repr-Kernels} can be used to construct
a multidimensional extension of the ns-form representation from Section~\ref{sec:Modified-ns-form},
using only one-dimensional quantities and norm estimates. This makes
our approach viable for $d>1$. We use the modified ns-form as in
the one-dimensional case described in Section~\ref{sec:Modified-ns-form}.
We find the ns-form essential for adaptive algorithms in more than
one dimension, since:

\begin{enumerate}
\item Scales do not interact as the operator is applied. All interactions
between scales are accounted for by the (inexpensive) projection at
the final step of the algorithm. 
\item For the same reason, the subdivision of space at different scales
naturally maps into the supporting data structures. We note that one
of the main difficulties in developing adaptive algorithms is in organizing
computations with blocks of an adaptive decomposition of a function
from different scales but with a common boundary. The methods for
such computations are known as mortar elements methods. In our approach
this issue does not present any obstacle, as all relevant interactions
are naturally accounted for by the data structures. 
\end{enumerate}
The key feature that makes our approach efficient in dimensions $d\ge2$
is the separated structure of the modified ns-form. Namely, the blocks
of $\mathbf{T}_{j}-\uparrow(\mathbf{T}_{j-1})$ are of the form (for
$d=3$)\begin{eqnarray}
\mathbf{T}_{j}^{\ell}-\uparrow(\mathbf{T}_{j-1}^{\ell}) & = & \sum_{m=1}^{M}w_{m}F^{j;ml_{1}}F^{j;ml_{2}}F^{j;ml_{3}}\nonumber \\
 &  & -\uparrow\left(\sum_{m=1}^{M}w_{m}F^{j-1;ml_{1}}F^{j-1;ml_{2}}F^{j-1;ml_{3}}\right)\nonumber \\
 & = & \sum_{m=1}^{M}w_{m}F^{j;ml_{1}}F^{j;ml_{2}}F^{j;ml_{3}}\label{eq:sep_apply}\\
 &  & -\sum_{m=1}^{M}w_{m}\uparrow(F^{j-1;ml_{1}})\cdot\uparrow(F^{j-1;ml_{2}})\cdot\uparrow(F^{j-1;ml_{3}}).\nonumber \end{eqnarray}
As in the case $d=1$, the norm of the operator blocks of $\mathbf{T}_{j}^{\ell}-\uparrow(\mathbf{T}_{j-1}^{\ell})$
decays rapidly with $\|\ell\|$, $\ell=(l_{1},l_{2},l_{3})$, and
the rate of decay depends on the number of vanishing moments of the
basis \cite{BE-CO-RO:1991}. Moreover, we limit the range of shift
indices $\|\ell\|$ using only one-dimensional estimates of the differences
\begin{equation}
F_{ii^{'}}^{j;\, m,2l}-\uparrow F_{ii^{'}}^{j;\, m,l}\,\,\,\mbox{and}\,\,\, F_{ii^{'}}^{j;\, m,2l+1}-\uparrow F_{ii^{'}}^{j;\, m,l}\label{eq:diff}\end{equation}
of operator blocks computed via (\ref{KernelFactor}). The norms of
individual blocks $F_{ii^{'}}^{j;\, m,l}$ are illustrated in Figure~\ref{fig:poisson_kernel_norms}
(a), for scale $j=1.$

By selecting the number of vanishing moments for a given accuracy,
it is sufficient to use $\|\ell\|_{\infty}\le2$ in practical applications
that we have encountered. Also, not all terms in the Gaussian expansion
of an operator need to be included since, depending on the scale $j$,
their contribution may be negligible for a given accuracy, as shown
in Figure~\ref{fig:poisson_kernel_norms} (b). Below we detail how
we select terms of the Gaussian expansion on a given scale as well
as the significant blocks of $\mathbf{T}_{j}^{\ell}-\uparrow(\mathbf{T}_{j-1}^{\ell})$.
This procedure establishes the band structure of the operator. We
then project the \emph{banded} operator $\uparrow(\mathbf{T}_{j-1}^{b_{j}})$
back to the scale $j-1$ and then combine blocks on the natural scale
for each projection in order to apply the operator efficiently as
was explained in Section~\ref{sec:Modified-ns-form} for the one-dimensional
case.

We note that in deciding which terms to keep in (\ref{eq:sep_apply}),
we do not compute the difference between the full three dimensional
blocks as it would carry a high computational cost; instead we use
estimates based on the one dimensional blocks of the separated representation.
We note that since the resulting band structure depends only on the
operator and the desired accuracy of its approximation, one of the
options is to store such band information as it is likely to be reused. 

In order to efficiently identify the significant blocks in $\mathbf{T}_{j}^{\ell}-\uparrow(\mathbf{T}_{j-1}^{\ell})$
as a function of $\ell$, we develop norm estimates based only on
the one-dimensional blocks. The difference between two terms of the
separated representation, say $F_{1}F_{2}F_{3}-G_{1}G_{2}G_{3}$,
may be written as \[
F_{1}F_{2}F_{3}-G_{1}G_{2}G_{3}=(F_{1}-G_{1})F_{2}F_{3}+G_{1}(F_{2}-G_{2})F_{3}+G_{1}G_{2}(F_{3}-G_{3}).\]
We average six different combinations of the three terms to include
all directions in a symmetric manner, which results in the norm estimate
\begin{eqnarray}
\| F_{1}F_{2}F_{3}-G_{1}G_{2}G_{3}\| & \le & \frac{1}{6}\mbox{sym}\left[\| F_{1}-G_{1}\|\| F_{2}\|\| F_{3}\|+\right.\label{eq:norm_estimate_ndim}\\
 &  & \left.\| G_{1}\|\| F_{2}-G_{2}\|\| F_{3}\|+\| G_{1}\|\| G_{2}\|\| F_{3}-G_{3}\|\right],\nonumber \end{eqnarray}
where the symmetrization is over the three directions and generates
$18$ terms. For the rotationally symmetric operators with Gaussian
expansion as in (\ref{KernelApproxByGaussians}) computing the right
hand side in this estimate involves just three types of one dimensional
blocks and their norms, \begin{eqnarray}
N_{\mbox{dif}}^{j;m;l} & = & \left\Vert F^{j;m;l}-\uparrow(F^{j-1;m;l})\right\Vert ,\nonumber \\
N_{F}^{j;m;l} & = & \left\Vert F^{j;m;l}\right\Vert ,\label{BlockNorms}\\
N_{\uparrow F}^{j;m;l} & = & \left\Vert \uparrow(F^{j-1;m;l})\right\Vert ,\nonumber \end{eqnarray}
where index $j$ indicates the scale, $m$ the term in the Gaussian
expansion (\ref{KernelApproxByGaussians}), and $l$ the position
of the block in a given direction.

These estimates allow us to discard blocks whose norm falls below
a given threshold of accuracy, namely, for each multi-index $\ell$,
we estimate\begin{eqnarray}
\left\Vert \mathbf{T}_{j}^{\ell}-\uparrow(\mathbf{T}_{j-1}^{\ell})\right\Vert  & \leq & \frac{1}{6}\sum_{m=1}^{M}w_{m}\,\,\mbox{sym}\left[N_{\mbox{dif}}^{j;m;l_{1}}N_{F}^{j;m;l_{2}}N_{F}^{j;m;l_{3}}+\right.\label{NormEstimateUsing1DBlocks}\\
 &  & \left.N_{\uparrow F}^{j;m;l_{1}}N_{\mbox{dif}}^{j;m;l_{2}}N_{F}^{j;m;l_{3}}+N_{\uparrow F}^{j;m;l_{1}}N_{\uparrow F}^{j;m;l_{2}}N_{\mbox{dif}}^{j;m;l_{3}}\right].\nonumber \end{eqnarray}

For each scale $j$ and each block $\mathbf{T}_{j}^{\ell}-\uparrow(\mathbf{T}_{j-1}^{\ell})$
labeled by the multi-index $\ell=(l_{1},l_{2},l_{3})$, we compute
all terms of the sum in (\ref{NormEstimateUsing1DBlocks}) and identify
the range $[m_{1},m_{2}]$ which we need to keep for that block, by
discarding from the sum terms whose cumulative contribution is below
$\epsilon$. If the entire sum falls below $\epsilon$, this range
may be empty and the entire $\mathbf{T}_{j}^{\ell}-\uparrow(\mathbf{T}_{j-1}^{\ell})$
is discarded. The range $[m_{1},m_{2}]$ differs significantly depending
whether or not the block is affected by the singularity of the kernel
as is illustrated in Figure~\ref{fig:poisson_kernel_norms}. In Figure~\ref{fig:poisson_kernel_norms}~(a)
the rate of decay for the blocks with shift $|l|=2,3$ is significantly
faster than for the blocks with $|l|\le1$ affected by the singularity.
We note that all blocks of the first 150 terms in the separated representation
(\ref{KernelApproxByGaussians2}) have norm 1 (and rank 1 as matrices)
and are not shown in Figure~\ref{fig:poisson_kernel_norms}~(a). 

Since the difference in (\ref{NormEstimateUsing1DBlocks}) involves
blocks upsampled from a coarser scale, all shifts $|l|\le3$ are affected
by the singularity. Figure \ref{fig:poisson_kernel_norms}~(b) shows
the r.h.s. of the estimate in (\ref{NormEstimateUsing1DBlocks}) for
different shifts along one of the directions, where the blocks along
the other two directions are estimated by the maximum norm over all
possible shifts. 

After discarding blocks with norms less that $\epsilon$ using the
estimate in (\ref{NormEstimateUsing1DBlocks}), we downsample the
remaining blocks of $\uparrow(\mathbf{T}_{j-1}^{\ell})$ back to the
original scale. This leaves only blocks of $\mathbf{T}_{j}^{\ell}$
on the scale $j$ and we remove additional blocks of $\mathbf{T}_{j}^{\ell}$
for the shifts $|l|=2,3$ where the decay is sufficiently fast to
make their contribution less than $\epsilon$.

This leads us to arrange the blocks on each scale into several subsets
by the effect the singularity of the kernel has on them and find the
appropriate range $[m_{1},m_{2}]$ for each subset. There are three
such sets in dimension $d=2$ and four sets in dimension $d=3$. For
each index $l,$ we will say that the index belongs to the \emph{core}
if $l=-1,0,1$ and to the \emph{boundary} otherwise. The core indices
correspond to one-dimensional blocks whose defining integrals include
the singularity of the kernel. We then divide all possible values
of the multi-index $\ell=(l_{1},l_{2},l_{3})$, according to the number
of core indices it has. In $d=3$ this gives us four sets:

\begin{itemize}
\item Core: all indices $(l_{1},l_{2},l_{3})$ belong to the core. 
\item Boundary-1: two of the indices belong to the core and one to the boundary
\item Boundary-2: one of the indices belongs to the core, the other two
to the boundary
\item Boundary-3: all indices $(l_{1},l_{2},l_{3})$ belong to the boundary.
\end{itemize}
We then find the range $[m_{1},m_{2}]$ for each subset and apply
blocks of each subset separately, thus avoiding unnecessary computations
with blocks whose contribution is negligible. This range analysis
only needs to be done once per operator and the desired accuracy,
and the results may be saved for repeated use.

\section{\label{sec:Multidim-application-nsform}Multidimensional adaptive
application of ns-form}

In this section we present an algorithm for applying the modified
non-standard form which is an extension of (\ref{alg:op-apply-1d})
(based on (\ref{eq:g_scales_band}) and (\ref{eq:apply_algo})) to
higher dimensions. We are now seeking to compute\[
g(\mathbf{x})=[Tf](\mathbf{x})=\int K(\mathbf{y}-\mathbf{x})f(\mathbf{y})d\mathbf{y},\]
where $\mathbf{x},\mathbf{y}\in\mathbb{R}^{d}$ for $d=2,3$. The
separated approximation (\ref{eq:sep_apply}) reduces the complexity
of applying the operator by allowing partial factorization of the
nested loops in each scale indicated by the order of summation and
illustrated for $d=3$,\begin{eqnarray*}
g^{j;l_{1}l_{2}l_{3}} & = & \sum_{m}w_{m}\left\{ \sum_{l'_{1}}F^{j;m;l_{1}-l'_{1}}\sum_{l'_{2}}F^{j;m,l_{2}-l'_{2}}\sum_{l'_{3}}F^{j;m,l_{3}-l'_{3}}\right.-\\
 &  & \left.\uparrow\left[\sum_{l'_{1}}F^{j-1;m;l_{1}-l'_{1}}\sum_{l'_{2}}F^{j-1;m,l_{2}-l'_{2}}\sum_{l'_{3}}F^{j-1;m,l_{3}-l'_{3}}\right]\right\} f^{j;l'_{1}l'_{2}l'_{3}}.\end{eqnarray*}

As described in the previous section, this evaluation is done by regions
of indices. These regions of indices are organized so that (for a
given accuracy) the number of retained terms of the separated representation
is roughly the same for all blocks within each region. Thus, we perform
the summation over the terms of the separated representation last,
applying only the terms that actually contribute to the result above
the requested accuracy threshold, according to estimate (\ref{NormEstimateUsing1DBlocks}).
Therefore, we avoid introducing checks per individual block and the
resulting loss of performance.

Just as in the one-dimensional case, we use (\ref{eq:sep_apply})
in a `natural scale' manner. That is, blocks belonging to scale $j$
are only applied on that scale. As in one-dimensional case, the interaction
between scales is achieved by the projection (\ref{eq:apply_algo})
that redistributes blocks accumulated in this manner properly between
the scales to obtain the adaptive tree for the resulting function.
The overall approach is the same as described in (\ref{alg:op-apply-1d}).
We note that, as expected, the separated representation requires more
detailed bookkeeping when constructing the data structures for the
operator.

\begin{rem}
\emph{Our multiresolution decomposition corresponds to the geometrically
varying refinement in finite element methods. In this case the adjoining
boxes do not necessarily share common vertices, forming what corresponds
to the so-called non-conforming grid. In finite element methods such
situation requires additional construction provided by the mortar
element methods. Mortar element methods were introduced by Patera
and his associates, see e.g. \cite{MA-MA-PA:1989,A-M-M-P:1990,BE-MA-PA:1993,BE-MA-PA:1994}.
These methods permit coupling discretizations of different types in
non-overlapping domains. Such methods are fairly complicated and involve,
for example, the introduction of interface conditions through an $L^{2}$
minimization. In our approach we do not face these issues at all and
do not have to introduce any additional interface conditions. The
proper construction for adjoining boxes is taken care by the redundant
tree data structure and Step~2 of Algorithm~\ref{alg:op-apply-nd}
for applying the kernel, which generates the necessary missing boxes
on appropriate scales.}
\end{rem}

\begin{rem}
\emph{Although Algorithm~\ref{alg:op-apply-nd} applies convolution
operators, only minor changes are needed to use it for non-convolutions.
Of course in such case, the separated representation for the modified
ns-form should be constructed by a different approach.}
\end{rem}
\begin{algorithm}[H]

\caption{\label{alg:op-apply-nd}Adaptive non-standard form operator application
in multiple dimensions (illustrated for $d=2$), $g=\mathbf{T}f$}

\begin{algorithmic}

\STATE {\bf Initialization:} Construct the redundant tree for $f$ and copy it
as skeleton tree for $g$ (see Section~\ref{sub:Redundant-tree-structure}).

\FORALL{scales $j=0,\ldots,n-1$}

  \FORALL{function blocks $g_{l_1l_2}^j$ in the tree for $g$ at scale $j$}

  \STATE \COMMENT{Step 1. Determine the list of all Core, Boundary-1 and
  Boundary-2 contributing operator blocks of the modified ns-form
  $F^{j;m;l_1-l_1'},F^{j;m;l_2-l_2'}$ (see
  Section~\ref{sec:Modified-ns-form-Nd}):}

  \IF{$g_{l_1l_2}^j$ belongs to a \emph{leaf}}

    \STATE Read operator blocks $F^{j;m;l_1-l_1'},F^{j;m;l_2-l_2'}$ for Core,
    Boundary-1 and Boundary-2, and their weights $w_m$ and corresponding
    ranges from the separated representation.

  \ELSE

    \STATE Read operator blocks $F^{j;m;l_1-l_1'},F^{j;m;l_2-l_2'}$ for
    Boundary-1 and Boundary-2, and their weights $w_m$ and corresponding ranges
    from the separated representation.

  \ENDIF

  \STATE \COMMENT{Step 2. Find the required blocks ${f_{l_1'l_2'}^j}$ of the
  input function $f$:}

  \IF {function block ${f_{l_1'l_2'}^j}$ exists in the redundant tree for $f$}

    \STATE Retrieve it.

  \ELSE

    \STATE Create it by interpolating from a coarser scale and cache for
    reuse.

  \ENDIF

  \STATE \COMMENT{Step 3. Output function block computation:}

  \STATE For each set $S$ of indices determined in \em{Step~1} (Core,
  Boundary-1, Boundary-2) and the corresponding ranges of terms in the
  separated representation, compute the sum

  \[
        \hat{g}_{l_1l_2}^{j;S}=\sum_m w_m \sum_{l_1'} F^{j;m;l_1-l_1'}
        \sum_{l_2'} F^{j;m;l_2-l_2'} f_{l_1'l_2'}^j,
  \]

  \STATE Add all computed sums to obtain $\hat{g}_{l_1l_2}^j$.

  \ENDFOR

\ENDFOR

\STATE \COMMENT{Step 4. Adaptive projection:}

\STATE Project resulting output function blocks $\hat{g}_{l_1l_2}^j$ on all
scales into a proper adaptive tree by using Eq.~(\ref{eq:g_assemble}).

\STATE Discard from the resulting tree unnecessary function blocks at the
requested accuracy.

\STATE {\bf Return:} the function $g$ represented by its adaptive tree.

\end{algorithmic}

\end{algorithm}

\subsection{Operation count estimates}

\subsubsection*{Adaptive decomposition of functions}

The cost of adaptively decomposing a function in $d$ dimensions is
essentially that of an adaptive wavelet transform. Specifically, it
takes $\mathcal{O}(N_{\textrm{blocks}}\cdot p^{d+1})$ operations
to compute such representation, where $N_{\textrm{blocks}}$ is the
final number of significant blocks in the representation and $p$
is the order of multiwavelet basis chosen. In comparison with the
usual wavelet transform, it appears to be significantly more expensive.
However, these $\mathcal{O}(p^{d+1})$ operations process $\mathcal{O}(p^{d})$
points, thus in counting significant coefficients as it is done in
the usual adaptive wavelet transform, we end up with $\mathcal{O}(p)$
operations per point.

\subsubsection*{Operator application}

The cost of applying an operator in the modified ns-form is $\mathcal{O}(N_{\textrm{blocks}}Mp^{d+1})$,
where $N_{\textrm{blocks}}$ is the number of blocks in the adaptive
representation of the input function, $M$ is the separation rank
of the kernel in (\ref{KernelApproxByGaussians}) and $p$ is the
order of the multiwavelet basis. For a given desired accuracy $\epsilon$,
we typically select $p\propto\log\epsilon^{-1}$; $M$ has been shown
to be proportional to $(\log\epsilon^{-1})^{\nu}$, where $\nu$ depends
on the operator \cite{BEY-MON:2005}. In our numerical experiments,
$M$ is essentially proportional to $\log\epsilon^{-1}$, since we
never use the full separated representation, as discussed in Section~\ref{sec:Modified-ns-form-Nd}
and illustrated in Figure~\ref{fig:poisson_kernel_norms}.

This operation count can be potentially reduced to $\mathcal{O}(N_{\textrm{blocks}}Mp^{d})$
by using the structure of the matrices in \ref{KernelFactor}, and
we plan to address this in the future.

\subsubsection*{Final projection}

After the operator has been applied to a function in a scale-independent
fashion, a final projection step is required as discussed in Section~\ref{sec:Modified-ns-form}.
This step requires $\mathcal{O}(N_{\textrm{blocks}}\cdot p^{d})$
operations, the same as in the original function decomposition. In
practice, this time is negligible compared to the actual operator
application.

\section{\label{sec:Numerical-examples}Numerical examples}

\subsection{The Poisson equation}

We illustrate the performance of the algorithm by solving the Poisson
equation in $d=3$\begin{equation}
\nabla^{2}\phi(\mathbf{r})=-\rho(\mathbf{r})\label{PoissonEquationExample}\end{equation}
with free space boundary conditions, $\phi(\mathbf{r})\to0$ and $\partial\phi/\partial r\to0$
as $r\to\infty$. We write the solution as \[
\phi(\mathbf{r})=\frac{1}{4\pi}\int_{\mathbb{R}^{3}}\frac{1}{\left|\mathbf{r}-\mathbf{r}'\right|}\rho(\mathbf{r}')d\mathbf{r}'\]
and adaptively evaluate this integral. We note that our method can
equally be used for $d=2,$ since the corresponding Green's function
can also be represented as a sum of Gaussians, and the operator application
algorithm has been implemented in for both $d=2$ and $d=3$. 

For our test we select \[
\phi(\mathbf{r})=\sum_{i=1}^{3}e^{-\alpha\left|\mathbf{r}-\mathbf{r}'\right|^{2}},\]
so that we solve the Poisson equation with\[
\rho(\mathbf{r})=-\nabla^{2}\phi(\mathbf{r})=-\sum_{i=1}^{3}(4\alpha^{2}\left|\mathbf{r}-\mathbf{r}_{i}\right|^{2}-6\alpha)e^{-\alpha\left|\mathbf{r}-\mathbf{r}_{i}\right|^{2}}.\]

Our parameters are chosen as follows: $\alpha=300$, $r_{1}=(0.5,0.5,0.5),$
$r_{2}=(0.6,0.6,0.5)$ and $r_{3}=(0.35,0.6,0.5).$ These ensure that
$\rho(r)$ is well below our requested thresholds on the boundary
of the computational domain. All numerical experiments were performed
on a Pentium-4 running at 2.8~GHz, with 2~GB of RAM. The results
are summarized in Table~\ref{table:poisson_3d}.

In order to gauge the speed of algorithm in reasonably computer-independent
terms, we use a similar approach to that of \cite{ETH-GRE:2001} and
also provide timings of the Fast Fourier Transform (FFT). Specifically,
we display timings for two FFTs as an estimate of the time needed
to solve the Poisson equation with a smooth right hand side and periodic
boundary conditions in a cube. As in \cite{ETH-GRE:2001}, we compute
the rate that estimates the number of processed points per second.
We observe that for our adaptive algorithm such rate varies between
$3.4\cdot10^{4}$ and $1.1\cdot10^{5}$ (see Table~\ref{table:poisson_3d}),
whereas for the FFTs it is around $10^{6}$ (see Table~\ref{table:fft_3d}).
We note that our algorithm is not fully optimized, namely, we do not
use the structure of the matrices in (\ref{KernelFactor}) and the
symmetries afforded by the radial kernels. We expect a substantial
impact on the speed by introducing these improvements and will report
them separately. 

We note that the multigrid method (see e.g. \cite{BRANDT:1977,BRANDT:1991})
is frequently used as a tool for solving the Poisson equation (and
similar problems) in differential form. The FFT-based gauge suggested
in \cite{ETH-GRE:2001} is useful for comparisons with these algorithms
as well.

\begin{table}
\begin{centering}Requested $\epsilon=10^{-3}$\par\end{centering}

\begin{centering}\begin{tabular}{|c|c|c|c|c|}
\hline 
$p$&
$E_{2}$&
$E_{\infty}$&
Time (s)&
Rate (pts/s)\tabularnewline
\hline
\hline 
6&
$5.0\cdot10^{-3}$&
$7.9\cdot10^{-1}$&
$1.2$&
$7.2\cdot10^{4}$\tabularnewline
\hline 
8&
$1.7\cdot10^{-3}$&
$1.2\cdot10^{-1}$&
$0.51$&
$7.3\cdot10^{4}$\tabularnewline
\hline 
10&
$4.4\cdot10^{-4}$&
$3.7\cdot10^{-2}$&
$0.68$&
$1.1\cdot10^{5}$\tabularnewline
\hline
\end{tabular}\par\end{centering}

\bigskip{}

\begin{centering}Requested $\epsilon=10^{-6}$\par\end{centering}

\begin{centering}\begin{tabular}{|c|c|c|c|c|}
\hline 
$p$&
$E_{2}$&
$E_{\infty}$&
Time (s)&
Rate (pts/s)\tabularnewline
\hline
\hline 
10&
$4.7\cdot10^{-6}$&
$3.6\cdot10^{-4}$&
10.3&
$5.7\cdot10^{4}$\tabularnewline
\hline 
12&
$8.5\cdot10^{-6}$&
$4.3\cdot10^{-5}$&
13.5&
$7.5\cdot10^{4}$\tabularnewline
\hline 
14&
$6.9\cdot10^{-8}$&
$5.2\cdot10^{-6}$&
20.0&
$8.0\cdot10^{4}$\tabularnewline
\hline
\end{tabular}\par\end{centering}

\bigskip{}

\begin{centering}Requested $\epsilon=10^{-9}$\par\end{centering}

\begin{centering}\begin{tabular}{|c|c|c|c|c|}
\hline 
$p$&
$E_{2}$&
$E_{\infty}$&
Time (s)&
Rate (pts/s)\tabularnewline
\hline
\hline 
16&
$2.5\cdot10^{-10}$&
$2.2\cdot10^{-8}$&
68.1&
$3.5\cdot10^{4}$\tabularnewline
\hline 
18&
$7.7\cdot10^{-11}$&
$3.5\cdot10^{-9}$&
100.3&
$3.4\cdot10^{4}$\tabularnewline
\hline 
20&
$1.2\cdot10^{-10}$&
$1.8\cdot10^{-8}$&
133.4&
$3.5\cdot10^{4}$\tabularnewline
\hline
\end{tabular}\par\end{centering}

\bigskip{}

\caption{\label{table:poisson_3d}Accuracy and timings for the adaptive solution
of the Poisson equation in (\ref{PoissonEquationExample}).}
\end{table}

\begin{table}
\begin{centering}\begin{tabular}{|c|c|c|}
\hline 
size&
Time (s)&
Rate (pts/s)\tabularnewline
\hline
\hline 
$32^{3}$&
0.02&
$1.7\cdot10^{6}$\tabularnewline
\hline 
$64^{3}$&
0.22&
$1.2\cdot10^{6}$\tabularnewline
\hline 
$128^{3}$&
2.21&
$9.5\cdot10^{5}$\tabularnewline
\hline 
$256^{3}$&
20.7&
$8.1\cdot10^{5}$\tabularnewline
\hline
\end{tabular}\par\end{centering}

\bigskip{}

\caption{\label{table:fft_3d}Timings of two 3D FFTs to estimate the speed
of a non-adaptive, periodic Poisson solver on a cube for smooth functions. }
\end{table}

\begin{figure}
\begin{centering}

  \includegraphics[width=0.35\textwidth]{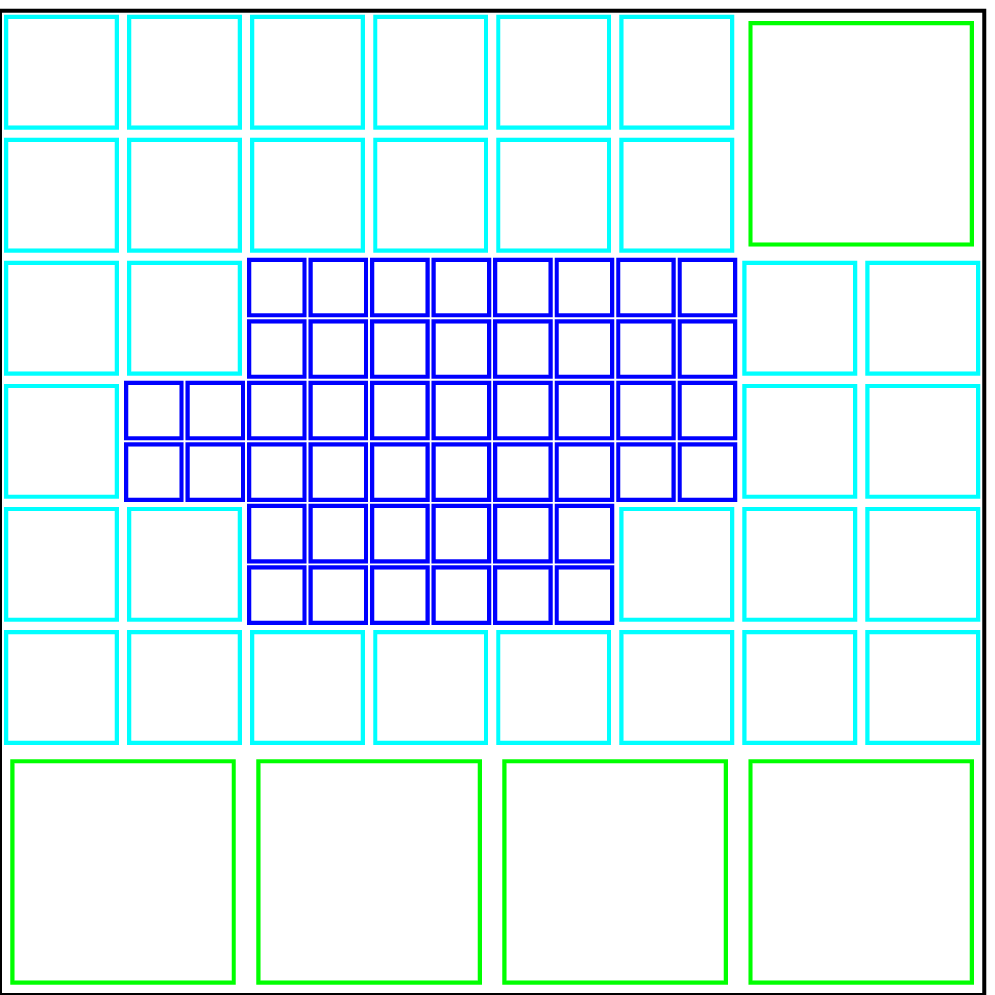} \hspace{5mm}%
  \includegraphics[width=0.6\textwidth]{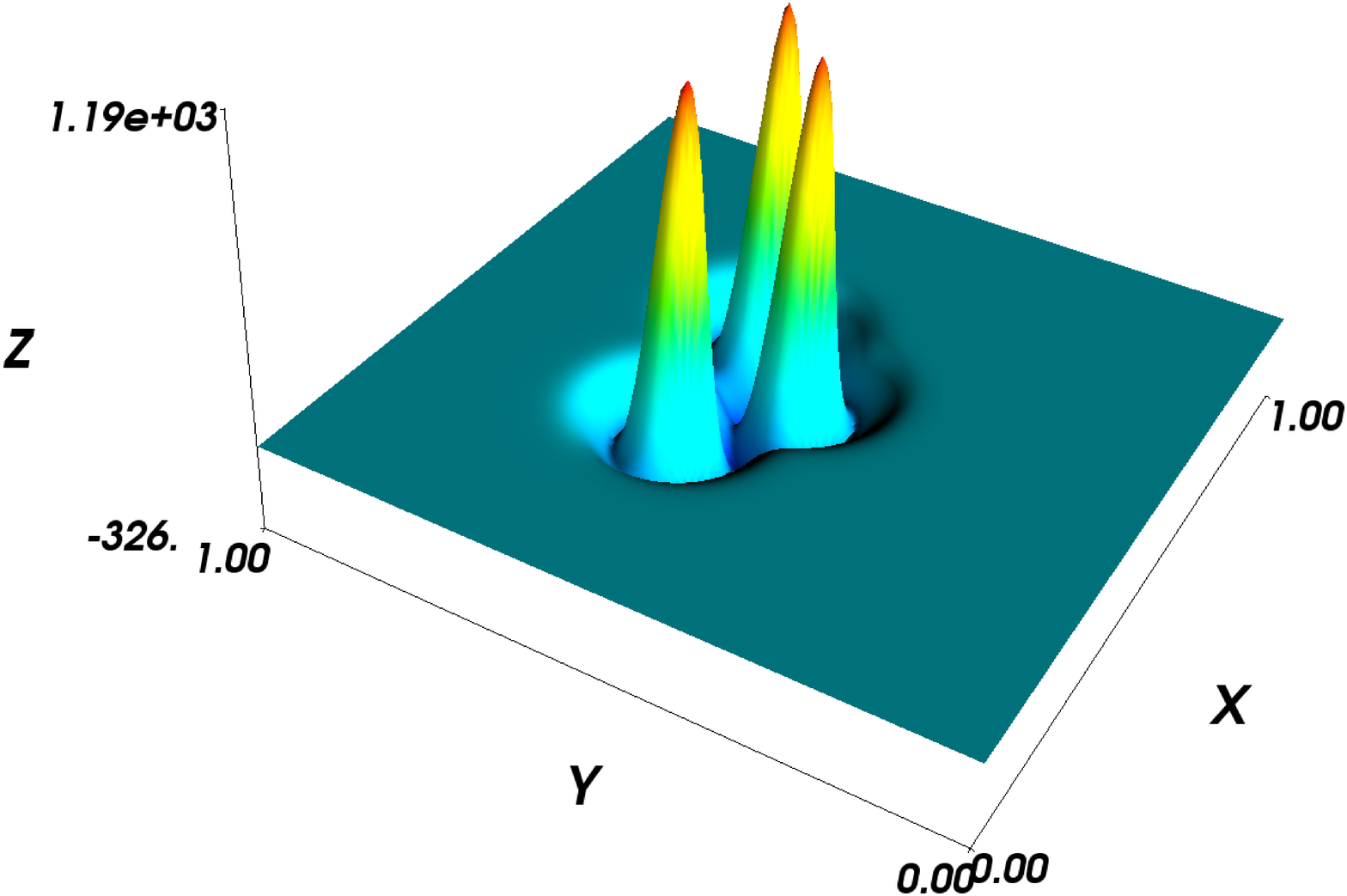}%
  \par
\end{centering}

\caption{\label{fig:poisson_3d}A two-dimensional slice of the three-dimensional
subdivision of space by the scaling functions and an illustration
of the source term for the Poisson equation (\ref{PoissonEquationExample}).}
\end{figure}

\subsection{The ground state of the hydrogen atom}

\begin{figure}

\begin{centering}

  \includegraphics[width=0.6\textwidth]{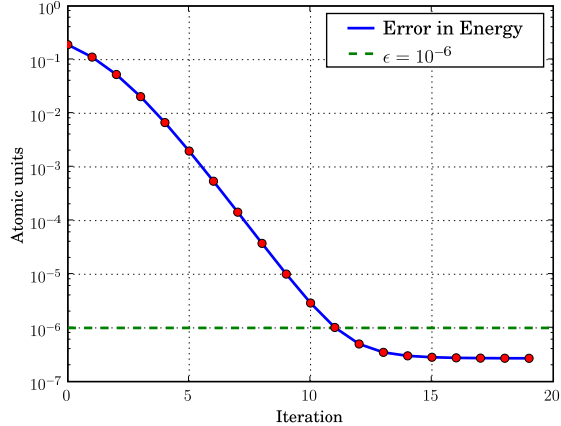}

  \par

\end{centering}

\caption{\label{fig:hydrogen_energy_error} Convergence of the iteration to
obtain the ground state of the hydrogen atom computed via formulation
in (\ref{Kalos}) for the non-relativistic Schr\"{o}dinger equation.
The requested accuracy in applying the Green's function is set to
$10^{-6}$. }
\end{figure}

A simple example of computing the ground state of the hydrogen atom
illustrates the numerical performance of the algorithms developed
in this paper, and their utility for constructing more complex codes.
The eigenfunctions $\psi$ for the hydrogen atom satisfy the time-independent
Schr\"{o}dinger equation (written in atomic units and spherical coordinates),\begin{equation}
-\frac{1}{2}\Delta\psi-\frac{1}{r}\psi=E\psi,\label{Schrodinger}\end{equation}
where $E$ is the energy eigenvalue. For the ground state, $E=-1/2$
and the (unnormalized) wave function is $\psi=e^{-r}$. Following
\cite{KALOS:1962}, we write \begin{equation}
\phi=-2G_{\mu}V\phi,\label{Kalos}\end{equation}
where $G_{\mu}=(-\Delta+\mu^{2}\mathcal{I})^{-1}$ is the Green's
function for some $\mu$ and $V=-1/r$ is the nuclear potential. For
$\mu=\sqrt{-2E}$ the solution $\phi$ of (\ref{Kalos}) has $\|\phi\|_{2}=1$
and coincides with that of (\ref{Schrodinger}). We solve (\ref{Kalos})
by a simple iteration starting from some value $\mu_{0}$ and changing
$\mu$ to obtain the solution with $\|\phi\|_{2}=1$. The algorithm
proceeds as follows:

\begin{enumerate}
\item Initialize with some value $\mu_{0}$ and function $\phi$. The number
of iterations of the algorithm is only weakly sensitive to these choices.
\item Compute the product of the potential $V$ and the function $\phi$.
\item Apply the Green's function $G_{\mu}$ to the product $V\phi$ via
the algorithm of this paper to compute \[
\phi_{new}=-2G_{\mu}V\phi.\]

\item Compute the energy for $\phi_{new}$,\[
E_{new}=\frac{\frac{1}{2}\langle\nabla\phi_{new},\nabla\phi_{new}\rangle+\langle V\phi_{new},\phi_{new}\rangle}{\langle\phi_{new},\phi_{new}\rangle}.\]

\item Set $\mu=\sqrt{-2E_{new}}$, $\phi=\phi_{new}/\|\phi_{new}\|$ and
return to Step~2.
\end{enumerate}
The iteration is terminated as the change in $\mu$ and $\|\phi\|-1$
falls below the desired accuracy. The progress of the iteration is
illustrated in Figure~\ref{fig:hydrogen_energy_error}. The computations
in Steps~2 and 4 use the three dimensional extension of the approach
described in \cite{A-B-G-V:2002}, to compute point-wise multiplications
of adaptively represented functions and weak differential operators
of the same.

This example illustrates an application of our algorithm to problems
in quantum chemistry. Multiresolution quantum chemistry developed
in \cite{H-F-Y-B:2003,H-F-Y-G-B:2004,Y-F-G-H-B:2004a,Y-F-G-H-B:2004}
also uses separated representations. The main technical difference
with \cite{H-F-Y-B:2003,H-F-Y-G-B:2004,Y-F-G-H-B:2004a,Y-F-G-H-B:2004}
is that we use the modified non-standard form and apply operator blocks
on their {}``natural'' scale, thus producing a fully adaptive algorithm.
We are currently using this algorithm as part of a new method for
solving the multiparticle Schr\"{o}dinger equation and will report
the results elsewhere.

\section{\label{sec:Discussion-and-conclusions}Discussion and conclusions}

We have shown that a combination of separated and multiresolution
representations of operators yields a new multidimensional algorithm
for applying a class of integral operators with radial kernels. We
note that the same algorithm is used for all such operators as they
are approximated by a weighted sum of Gaussians. This fact makes our
approach applicable across multiple fields, where a single implementation
of the core algorithm can be reused for different specific problems.
The algorithm is fully adaptive and avoids issues usually addressed
by mortar methods. 

The method of approximation underlying our approach is distinct from
that of the FMM, has similar efficiency and has the advantage of being
more readily extendable to higher dimensions. We also note that semi-analytic
approximations via weighted sums of Gaussians provide additional advantages
in some applications. Although we described the application of kernels
in free space, there is a simple extension to problems with radial
kernels subject to periodic, Dirichlet or Neumann boundary conditions
on a cube that we will describe separately.

The algorithm may be extended to classes of non-convolution operators,
e.g., the Calderon-Zygmund operators. For such extensions the separated
representation may not be available in analytic form, as it is for
the operators of this paper, and may require a numerical construction.
The separated representation of the kernel permits further generalization
of our approach to dimensions $d\gg3$ for applying operators to functions
in separated representation.

A notable remaining challenge is an efficient, high order extension
of this technique to the application of operators on domains with
complicated geometries and surfaces.

\section{\label{sec:Appendix}Appendix }

\subsection{Scaling functions}

We use either the Legendre polynomials $P_{0},\ldots,P_{p-1}$ or
the interpolating polynomials on the Gauss-Legendre nodes in $[-1,1]$
to construct an orthonormal basis for each subspace $\mathbf{V}_{j}$
\cite{ALPERT:1993,A-B-G-V:2002}. 

Let us briefly describe some properties of the Legendre scaling functions
$\phi_{k}$, $k=0,\ldots,p-1$, defined as\begin{equation}
\phi_{k}(x)=\left\{ \begin{array}{ll}
\sqrt{2k+1}P_{k}(2x-1), & x\in[0,1]\\
0, & x\notin[0,1]\end{array}\right.,\label{scalingFNC}\end{equation}
and forming a basis for $\mathbf{V}_{0}.$ The subspace $\mathbf{V}_{j}$
is spanned by $2^{j}p$ functions obtained from $\phi_{0},\ldots,\phi_{p-1}$
by dilation and translation,\begin{equation}
\phi_{kl}^{j}(x)=2^{j/2}\phi_{k}(2^{j}x-l),\quad k=0,\ldots,p-1,\quad l=0,\ldots,2^{j}-1.\end{equation}
These functions have support on $[2^{-j}l,2^{-j}(l+1)]$ and satisfy
the orthonormality condition \begin{equation}
\int_{-\infty}^{\infty}\phi_{kl}^{j}(x)\phi_{k'l'}^{j}(x)dx=\delta_{kk'}\delta_{ll'}\,.\end{equation}
A function $f$, defined on $[0,1]$, is represented in the subspace
$\mathbf{V}_{j}$ by its normalized Legendre expansion\begin{equation}
f(x)=\sum_{l=0}^{2^{j}-1}\sum_{k=0}^{p-1}s_{kl}^{j}\phi_{kl}^{j}(x),\label{eq:func_interp_1d}\end{equation}
where the coefficients $s_{kl}^{j}$ are computed via\begin{equation}
s_{kl}^{j}=\int_{2^{-j}l}^{2^{-j}(l+1)}f(x)\phi_{kl}^{j}(x)\, dx.\end{equation}
As long as $f(x)$ is smooth enough and is well approximated on $[2^{-j}l,2^{-j}(l+1)]$
by a polynomial of order up to $2p-1$, we may use Gauss-Legendre
quadratures to calculate the $s_{kl}^{j}$ via\begin{equation}
s_{kl}^{j}=2^{-j/2}\sum_{i=0}^{p-1}f(2^{-j}(x_{i}+l))\phi_{k}(x_{i})w_{i},\label{eq:fvals2coefs}\end{equation}
 where $x_{0},\ldots,x_{p-1}$ are the roots of $P_{p}(2x-1)$ and
$w_{0,}\ldots,w_{p-1}$ are the corresponding quadrature weights.

In more than one dimension, the above formulas are extended by using
a tensor product basis in each subspace. For example, in two dimensions
equation~(\ref{eq:func_interp_1d}) becomes\begin{equation}
f(x,x')=\sum_{l=0}^{2^{j}-1}\sum_{k=0}^{p-1}\sum_{l'=0}^{2^{j}-1}\sum_{k'=0}^{p-1}s_{kk'll'}^{j}\phi_{kl}^{j}(x)\phi_{k'l'}^{j}(x').\label{eq:func_interp_2d}\end{equation}

\subsection{Multiwavelets}

We use piecewise polynomial functions $\psi_{0,}\ldots,\psi_{p-1}$
as an orthonormal basis for $\mathbf{W}_{0}$ \cite{ALPERT:1993,A-B-G-V:2002},\begin{equation}
\int_{0}^{1}\psi_{i}(x)\psi_{j}(x)dx=\delta_{ij}.\end{equation}
Since $\mathbf{W}_{0}\perp\mathbf{V}_{0}$, the first $p$ moments
of all $\psi_{0,}\ldots,\psi_{p-1}$ vanish:\begin{equation}
\int_{0}^{1}\psi_{i}(x)x^{i}dx=0,\quad i,j=0,1,\ldots,p-1.\end{equation}
The space $\mathbf{W}_{j}$ is spanned by $2^{j}p$ functions obtained
from $\psi_{0,}\ldots,\psi_{p-1}$ by dilation and translation,\begin{equation}
\psi_{kl}^{j}(x)=2^{j/2}\psi_{k}(2^{j}x-l),\quad k=0,\ldots,p-1,\quad l=0,\ldots,2^{j}-1,\end{equation}
and supported in the interval $I_{jl}=[2^{-j}l,2^{-j}(l+1)].$ A function
$f(x)$ defined on $[0,1]$ is represented in the multiwavelet basis
on $n$ scales by\begin{equation}
f(x)=\sum_{k=0}^{p-1}s_{k,0}^{0}\phi_{k}(x)+\sum_{j=0}^{n-1}\sum_{l=0}^{2^{j}-1}\sum_{k=0}^{p-1}d_{kl}^{j}\psi_{kl}^{j}(x)\label{eq:multiwavelet_expansion}\end{equation}
with the coefficients $d_{kl}^{j}$ computed via\begin{equation}
d_{kl}^{j}=\int_{2^{-j}l}^{2^{-j}(l+1)}f(x)\psi_{kl}^{j}(x)\, dx.\end{equation}

\subsection{Two-scale relations}

The relation between subspaces, $\mathbf{V}_{0}\oplus\mathbf{W}_{0}=\mathbf{V}_{1}$,
is expressed via the two-scale difference equations,\begin{eqnarray}
\phi_{k}(x) & = & \sqrt{2}\sum_{k'=0}^{p-1}\left(h_{kk'}^{(0)}\phi_{k'}(2x)+h_{kk'}^{(1)}\phi_{k'}(2x-1)\right),\quad k=0,\ldots,p-1,\\
\psi_{k}(x) & = & \sqrt{2}\sum_{k'=0}^{p-1}\left(g_{kk'}^{(0)}\phi_{k'}(2x)+g_{kk'}^{(1)}\phi_{k'}(2x-1)\right),\quad k=0,\ldots,p-1,\end{eqnarray}
where the coefficients $h_{ij}^{(0)}$, $h_{ij}^{(1)}$ and $g_{ij}^{(0)}$,
$g_{ij}^{(1)}$ depend on the type of polynomial basis used (Legendre
or interpolating) and its order $p$. The matrices of coefficients
\[
H^{(0)}=\{ h_{kk'}^{(0)}\},\quad H^{(1)}=\{ h_{kk'}^{(1)}\},\quad G^{(0)}=\{ g_{kk'}^{(0)}\},\quad G^{(1)}=\{ g_{kk'}^{(1)}\}\]
are the multiwavelet analogues of the quadrature mirror filters in
the usual wavelet construction, e.g., \cite{DAUBEC:1992}. These matrices
satisfy a number of important orthogonality relations and we refer
to \cite{A-B-G-V:2002} for complete details, including the construction
of the $H,\, G$ matrices themselves. Let us only state how these
matrices are used to connect the scaling $s_{kl}^{j}$ and wavelet
$d_{kl}^{j}$ coefficients on neighboring scales $j$ and $j+1$.
The \emph{decomposition procedure} ($j+1\rightarrow j$) is based
on \begin{eqnarray}
s_{kl}^{j} & = & \sum_{k'=0}^{p-1}\left(h_{kk'}^{(0)}s_{k,2l}^{j+1}+h_{kk'}^{(1)}s_{k,2l+1}^{j+1}\right),\label{eq:decomp_scaling}\\
d_{kl}^{j} & = & \sum_{k'=0}^{p-1}\left(g_{kk'}^{(0)}s_{k,2l}^{j+1}+g_{kk'}^{(1)}s_{k,2l+1}^{j+1}\right);\label{eq:decomp_wavelet}\end{eqnarray}
the \emph{reconstruction} ($j\rightarrow j+1$) is based on\begin{eqnarray}
s_{k,2l}^{j+1} & = & \sum_{k'=0}^{p-1}\left(h_{kk'}^{(0)}s_{kl}^{j}+g_{kk'}^{(0)}d_{kl}^{j}\right),\label{eq:recon_even}\\
s_{k,2l+1}^{j+1} & = & \sum_{k'=0}^{p-1}\left(h_{kk'}^{(1)}s_{kl}^{j}+g_{kk'}^{(1)}d_{kl}^{j}\right).\label{eq:recon_odd}\end{eqnarray}

\subsection{\label{sub:Cross-correlation-functions-of}Cross-correlation of the
scaling functions}

For convolution operators we only need to compute integrals with the
cross-correlation functions of the scaling functions, \begin{equation}
\Phi_{ii'}(x)=\int_{-\infty}^{\infty}\phi_{i}(x+y)\phi_{i'}(y)dy.\label{CrossCor}\end{equation}
 Since the support of the scaling functions is restricted to $[0,1]$,
the functions $\Phi_{ii'}$ are zero outside the interval $[-1,1]$
and are polynomials on $[-1,0]$ and $[0,1]$ of degree $i+i'+1$,

\begin{equation}
\Phi_{ii'}(x)\,=\left\{ \begin{array}{cr}
\Phi_{ii'}^{+}(x), & \,\,\,0\leq x\leq1,\\
\Phi_{ii'}^{-}(x), & -1\leq x<0,\\
0, & 1<|x|,\end{array}\right.\label{crosscor.00}\end{equation}
 where $i,i'=0,\ldots,p-1$ and\begin{equation}
\Phi_{ii'}^{+}(x)\,=\,\int_{0}^{1-x}\,\phi_{i}(x+y)\,\phi_{i'}(y)dy\,,\,\,\,\,\,\,\,\Phi_{ii'}^{-}(x)\,=\,\int_{-x}^{1}\,\phi_{i}(x+y)\,\phi_{i'}(y)dy\,.\label{crosscor.01}\end{equation}
 We summarize relevant properties of the cross-correlation functions
$\Phi_{ii'}$ in

\begin{prop}
\label{pro:CrossCor}~
\begin{enumerate}
\item Transposition of indices: $\Phi_{ii'}(x)=(-1)^{i+i'}\Phi_{i'i}(x)$,\\

\item Relations between $\Phi^{+}$and $\Phi^{-}$: $\Phi_{i,i'}^{-}(-x)=(-1)^{i+i'}\Phi_{i,i'}^{+}(x)$
for $0\leq x\leq1$,\\

\item Values at zero: $\Phi_{ii'}(0)=0$ if $i\ne i'$, and $\Phi_{ii}(0)=1$
for $i=0,\dots,p-1$,\\

\item Upper bound: $\max_{x\in[-1,1]}|\Phi_{ii'}(x)|\leq1$ for $i,i'=0,\dots,p-1$,\\

\item Connection with the Gegenbauer polynomials: \\
$\Phi_{00}^{+}(x)=\frac{1}{2}\, C_{1}^{(-1/2)}(2x-1)+\frac{1}{2}$
and $\Phi_{l0}^{+}(x)=\frac{1}{2}\,\sqrt{2l+1}\, C_{l+1}^{(-1/2)}(2x-1)$,
for $l=1,2,\dots$, where $C_{l+1}^{(-1/2)}$ is the Gegenbauer polynomial,\\

\item Linear expansion: if $i'\ge i$ then we have\begin{equation}
\Phi_{ii'}^{+}(x)=\sum_{l=i'-i}^{i'+i}\, c_{ii'}^{l}\Phi_{l0}^{+}(x),\label{expans.00}\end{equation}
 where\begin{equation}
c_{ii'}^{l}=\left\{ \begin{array}{lcc}
4l(l+1)\int_{0}^{1}\Phi_{ii'}^{+}(x)\Phi_{l0}^{+}(x)(1-(2x-1)^{2})^{-1}dx, &  & i'>i,\\
\\4l(l+1)\int_{0}^{1}(\Phi_{ii}^{+}(x)-\Phi_{00}^{+}(x))\,\Phi_{l0}^{+}(x)(1-(2x-1)^{2})^{-1}dx, & \mbox{} & i'=i,\end{array}\right.\label{CoefForPhi}\end{equation}
for $l\ge1$ and $c_{ii'}^{0}=\delta_{ii'}$.\\

\item Vanishing moments: we have $\int_{-1}^{1}\Phi_{00}(x)\, dx=1$ and
$\int_{-1}^{1}x^{k}\Phi_{ii'}(x)\, dx=0$ for $i+i'\geq1$ and $0\leq k\leq i+i'-1$. 
\end{enumerate}
\end{prop}
Proof of these properties can be found in \cite{B-C-F-H:2007}.

\subsection{\label{sub:Separated-representations}Separated representations of
radial functions}

As an example, consider approximating the function $1/r^{\alpha}$
by a collection of Gaussians. The number of terms needed for this
purpose is mercifully small. We have \cite{BEY-MON:2005}

\begin{prop}
\label{pro:trapez}For any $\alpha>0$, $0<\delta\leq1$, and $0<\epsilon\leq\min\left\{ \frac{1}{2},\frac{8}{\alpha}\right\} $,
there exist positive numbers $\tau_{m}$ and $w_{m}$ such that\begin{equation}
{\Big|}r^{-\alpha}-\sum_{m=1}^{M}w_{m}e^{-\tau_{m}r^{2}}{\Big|}\leq r^{-\alpha}\epsilon,\textrm{\mbox{ for all }$\delta\leq r\leq1$}\label{trapez.01}\end{equation}
 with\begin{equation}
M=\log\epsilon^{-1}[c_{0}+c_{1}\log\epsilon^{-1}+c_{2}\log\delta^{-1}],\label{trapez.01a}\end{equation}
 where $c_{k}$ are constants that only depend on $\alpha$. For fixed
power $\alpha$ and accuracy $\epsilon$, we have $M=\mathcal{O}(\log\delta^{-1})$. 
\end{prop}
The proof of Proposition~\ref{pro:trapez} in \cite{BEY-MON:2005}
is based on using the trapezoidal rule to discretize an integral representation
of $1/r^{\alpha}$. Similar estimates may be obtained for more general
radial kernels using their integral representations as in (\ref{IntegralRepr}).

We note that approximations of the function $1/r$ via sums of Gaussians
have been also considered in \cite{KUTZEL:1994,BRAESS:1995,BRA-HAC:2005}.

\subsection{Estimates}

By selecting appropriate $\epsilon$ and $\delta$ in the separated
representations of a radial kernel $K$ (as in Proposition~\ref{pro:trapez}),
we obtain a separated approximation for the coefficients\begin{equation}
t_{ii',jj',kk'}^{j;\,{\bf \ell}}\,=\,2^{-3j}\,\int_{[-1,1]^{3}}\, K(2^{-j}(\mathbf{x}+{\bf \ell}))\,\Phi_{ii'}(x_{1})\,\Phi_{jj'}(x_{2})\,\Phi_{kk'}(x_{3})\, d\mathbf{x}\,.\label{arbitrary-kernel}\end{equation}
Since the number of terms, $M$, depends logarithmically on $\epsilon$
and $\delta$, we achieve any finite accuracy with a very reasonable
number of terms. For example, for the Poisson kernel $K(r)=1/4\pi r$,
we have the following estimate \cite{B-C-F-H:2007},

\begin{prop}
\label{estimates} For any $\epsilon>0$ and $0<\delta\le1$ the coefficients
$t_{ii',jj',kk'}^{j;\,{\bf \ell}}$ in (\ref{arbitrary-kernel}) have
an approximation with a low separation rank, \begin{equation}
r_{ii',jj',kk'}^{j;\,{\bf \ell}}\,=\sum_{m=1}^{M}\, w_{m}\, F_{ii'}^{m,l_{1}}\, F_{jj'}^{m,l_{2}}\, F_{kk'}^{m,l_{3}}\,,\label{kern3d.0101}\end{equation}
 such that if $\max_{i}|l_{i}|\geq2$, then \begin{equation}
|t_{ii',jj',kk'}^{j;\,{\bf \ell}}\,-\, r_{ii',jj',kk'}^{j;\,{\bf \ell}}|\,\leq\, c_{0}2^{-2j}\epsilon,\label{away-from-singularity}\end{equation}
 and if $\max_{i}|l_{i}|\leq1$, then \begin{equation}
|t_{ii',jj',kk'}^{j;\,{\bf \ell}}\,-\, r_{ii',jj',kk'}^{j;\,{\bf \ell}}|\,\leq\,2^{-2j}(c_{1}\delta^{2}+\, c_{0}\epsilon)\label{near-singularity}\end{equation}
 where $\epsilon$, $\delta$, $M$, $\tau_{m}$, $w_{m}$, $m=1,\ldots,M$
are described in Proposition~\ref{pro:trapez} for $\alpha=1$ and
$c_{0}$ and $c_{1}$ are (small) constants.
\end{prop}
As described in Section~\ref{sec:Modified-ns-form-Nd}, our adaptive
algorithm selects only some of the terms, as needed on a given scale
for the desired accuracy $\epsilon$.

\subsection{Evaluation of integrals with the cross-correlation functions }

We need to compute integrals in (\ref{KernelFactor}), where the cross-correlation
functions $\Phi_{ii'}$ are given in (\ref{CrossCor}). We note that
using (\ref{expans.00}), it is sufficient to compute \[
F_{i0}^{j;m,l}=\frac{1}{2^{j}}\int_{-1}^{1}e^{-\tau_{m}(x+l)^{2}/4^{j}}\Phi_{i0}(x)dx\]
for $0\le i\le2p-1$ rather than $p^{2}$ integrals in (\ref{KernelFactor}).
Using the relation between $\Phi^{-}$ and $\Phi^{+}$ in Proposition~\ref{pro:CrossCor},
we have\[
\int_{-1}^{0}\Phi_{i0}^{-}(x)e^{-\tau_{m}(x+l)^{2}}dx=\int_{0}^{1}\Phi_{i0}^{-}(-x)e^{-\tau_{m}(-x+l)^{2}}dx=(-1)^{i}\int_{0}^{1}\Phi_{i0}^{+}(x)e^{-\tau_{m}(x-l)^{2}}dx,\]
so that \[
F_{i0}^{j;m,l}=\int_{0}^{1}[e^{-\tau_{m}(x+l)^{2}}+(-1)^{i}e^{-\tau_{m}(x-l)^{2}}]\Phi_{i0}^{+}(x)dx.\]


\end{document}